\documentclass[12pt]{article}
\usepackage[cp866]{inputenc}
\usepackage{citehack}
\usepackage{amsmath}
\usepackage{amsfonts}
\usepackage{amssymb}
\usepackage{longtable}
\newtheorem{theorem}{Theorem}[section]
\newtheorem{prop}{Proposition}[section]
\newtheorem{remark}{Remark}[section]
\newtheorem{definit}{Definition}[section]
\newtheorem{lem}{Lemma}[section]
\newtheorem{ex}{Example}[section]
\newtheorem{tab}{Table}
\newtheorem{corol}{Corollary}[section]

\def\spa{\mathop\text{{\rm span}}\nolimits}
\def\Hom{\mathop\text{\rm Hom}\nolimits}
\def\pr{\mathop\text{\rm pr}\nolimits}

\def\str{\mathop\text{\rm str}\nolimits}
\def\Ric{\mathop\text{\rm Ric}\nolimits}

\def\id{\mathop\text{\rm id}\nolimits}
\def\sign{\mathop\text{\rm sign}\nolimits}
\def\Ad{\mathop\text{\rm Ad}\nolimits}
\def\ad{\mathop\text{\rm ad}\nolimits}

\def\GL{\text{\rm GL}}

\def\Ort{\text{\rm O}}
\def\Un{\text{\rm U}}
\def\Sp{\text{\rm Sp}}
\def\SU{\text{\rm SU}}

\def\Real{\mathbb{R}}
\def\Co{\mathbb{C}}

\def\g{\mathfrak{g}}
\def\h{\mathfrak{h}}

\def\osp{\mathfrak{osp}}
\def\hosp{\mathfrak{hosp}}
\def\sl{\mathfrak{sl}}
\def\sp{\mathfrak{sp}}
\def\gl{\mathfrak{gl}}
\def\su{\mathfrak{su}}
\def\u{\mathfrak{u}}
\def\hol{\mathfrak{hol}}
\def\q{\mathfrak{q}}
\def\pe{\mathfrak{pe}}

\def\p{\partial}

\def\Ga{\Gamma}
\def\na{\nabla}

\def\t{\tilde}
\def\ti{\thicksim}

\def\al{\alpha}
\def\ga{\gamma}
\def\pa{{\partial_a}}
\def\pb{{\partial_b}}
\def\pai{{\partial_i}}

\def\pga{{\partial_\gamma}}

\def\T{{\cal T}}
\def\M{{\cal M}}
\def\N{{\cal N}}

\def\E{{\cal E}}
\def\F{{\cal F}}
\def\O{{\cal O}}
\def\Hols{{\mathcal Hol}}
\def\Hol{\mathop\text{\rm Hol}\nolimits}

\def\R{{\cal R}}

\def\bA{{\bar A}}
\def\bB{{\bar B}}
\def\bC{{\bar C}}
\def\bD{{\bar D}}
\def\hA{{\hat A}}
\def\hB{{\hat B}}
\def\hC{{\hat C}}
\def\hD{{\hat D}}

\title{Holonomy  of supermanifolds}

\author{Anton S. Galaev\footnote{Supported from the Basic Research Center no.
LC505 (Eduard \v{C}ech Center for Algebra and Geometry) of Ministry of Education, Youth and Sport of Czech Republic.} }

\begin{document}

\maketitle\vskip-50ex {\renewcommand{\abstractname}{Abstract}\begin{abstract} Holonomy groups and holonomy algebras for
connections on locally free sheaves over supermanifolds are introduced. A one-to-one correspondence between parallel
sections and holonomy-invariant vectors, and a one-to-one correspondence between parallel locally direct subsheaves  and
holonomy-invariant vector supersubspaces are obtained. As the special case, the holonomy of  linear connections on
supermanifolds is studied. Examples of parallel geometric structures on supermanifolds and the corresponding holonomies
are given. For Riemannian supermanifolds an analog of the Wu theorem is proved. Berger superalgebras are defined and their
examples are given.
\end{abstract}


{\bf Mathematical subject codes:} 58A50, 53C29

\tableofcontents

\section{Introduction}

The holonomy groups play a big role in the study of connections on vector bundles over smooth manifolds. They link
geometric and algebraic properties.  In particular, they allow  us to find parallel sections in geometric vector bundles
associated to the manifold, such as the tangent bundle, tensor bundles, or the spin bundle, as holonomy-invariant objects,
see \cite{Be,Bryant1,IRMA,Jo,Jo2007,K-N}.

In the present paper we introduce holonomy groups for connections on  supermanifolds.

In Section \ref{supermf}  some necessary preliminaries on supermanifolds are given. Section \ref{smoothhol} is an
introduction to the theory of holonomy for connections on vector bundles over smooth manifolds. In Section \ref{superhol}
we define the holonomy for connections on locally free sheaves over supermanifolds. First we define the holonomy algebra,
for this we generalize the Ambrose-Singer theorem and use the covariant derivatives of the curvature tensor and parallel
displacements. Then we define the holonomy group as a Lie supergroup.  In Section \ref{sectinf} we define the
infinitesimal holonomy algebra and  show that in the analytic settings it coincides with the holonomy algebra.
 In Section \ref{parallel} we
study parallel sections of sheaves over supermanifolds. It is shown that any parallel section is uniquely defined by its
value at any point, in spite of the fact that generally the sections of  sheaves over supermanifolds are not defined by
their values at all points. After this we obtain a one-to-one correspondence between parallel sections and
holonomy-invariant vectors, as in the case of vector bundles over smooth manifolds. The definition of the holonomy was
motivated by this correspondence. In Section \ref{sps} a one-to-one correspondence between parallel locally direct
subsheaves  and holonomy-invariant vector supersubspaces is obtained. Then we turn to study holonomy of linear connections
on supermanifolds. In Section \ref{LinConnect}  a one-to-one correspondence between parallel tensors on a supermanifold
and holonomy-invariant tensors at one point is obtained. We consider examples of parallel structures on supermanifolds and
give the equivalent conditions in terms of holonomy. In Section \ref{secBerger} Berger superalgebras are introduced. These
superalgebras generalize the usual Berger algebras and  they can be considered as candidates to holonomy algebras of
linear torsion-free connections on supermanifolds. In Section \ref{secsym}  holonomy of locally symmetric supermanifolds
is considered. In Section \ref{holRiem}  the case of the Levi-Civita connections on Riemannian supermanifolds is studied.
K\"ahlerian, special K\"ahlerian, hyper-K\"ahlerian and quaternionic-K\"ahlerian supermanifolds are characterized by their
holonomy. It is shown that special K\"ahlerian supermanifolds are Ricci-flat and, conversely, Ricci-flat simply connected
K\"ahlerian supermanifolds are special K\"ahlerian. A generalization of the Wu theorem is proved. In Section
\ref{exBerger} we give examples of complex Berger superalgebras.

Thus holonomy of supermanifolds that is introduced in the present paper is an appropriate generalization of the usual
holonomy of smooth manifolds, as  many properties are preserved.

\section{Supermanifolds}\label{supermf}

In this section we  give some necessary preliminaries on supermanifolds. An introduction to linear superalgebra and to the
theory of supermanifolds can be found in \cite{DelMor,Leites,Leites1,Manin,Var}.

{\it  A real smooth (analytic)  supermanifold} $\M$ of dimension $n|m$ is a pair $(M,\O_\M)$, where $M$ is a real smooth
(analytic) manifold of dimension $n$, $\O_\M$ is a sheaf of superalgebras over $\Real$ such that locally $$\O_\M(U)\simeq
\O_M(U)\otimes \Lambda_{\xi^1,...,\xi^m}.$$ Here $\O_M(U)$ is the algebra of smooth (analytic) functions on $U\subset M$
and $\Lambda_{\xi^1,...,\xi^m}$ is the Grassmann algebra of $m$ generators. If $m=0$, then $\M=M$.    The sections of the
sheaf $\O_\M$ are called {\it superfunctions (or just functions)} on $\M$. There exists  the canonical projection
$\ti:\O_\M\to \O_M$, $f\mapsto \t f$, where $\O_M$ is the sheaf of smooth functions on $M$. {\it The value of a
superfunction} $f$ at a point $x$ is $\t f(x)$.

We will use the following convention about the ranks of the indices $i,j,k=1,...,n$,  $\al,\beta,\ga=1,...,m$ and
$a,b,c=1,...,n+m$. We will use the Einstein rule for sums. Let $U\subset M$ be as above and let $x^i$ be local coordinates
on $U$, then the system $(U,x^i,\xi^\al)$ is called {\it a system of local coordinates} on $\M$ over $U$. We will denote
such system also by $(U,x^a)$, where $x^{n+\al}=\xi^\al$. For any $f\in\O_\M(U)$ we get
\begin{equation} f=\sum_{r=0}^m\sum_{\al_1<\cdots<\al_r}f_{\al_1...\al_r}\xi^{\al_1}\cdots\xi^{\al_r},\label{f}\end{equation} where
 $f_{\al_1...\al_r}\in\O_M(U)$ and $f_\emptyset=\t f$.  For any $\al_1<\cdots<\al_r$ and
 permutation  $\sigma:\{\al_1,...,\al_r\}\to \{\al_1,...,\al_r\}$ we assume that
 $f_{\sigma(\al_1)...\sigma(\al_r)}=\sign_\sigma f_{\al_1...\al_r}$. If two of the numbers $\al_1,...,\al_r$ are
 equal, we assume that $f_{\al_1...\al_r}=0$.

Denote by $\T_\M$ the tangent sheaf, i.e. the sheaf of superderivatives of the sheaf $\O_\M$. If $(U,x^i,\xi^\al)$ is a
system of local coordinates, then the vector fields $\p_{x^i},\p_{\xi^\al}$ form a basis of the supermodule $\T_\M(U)$
over the superalgebra $\O_\M(U)$. The vector fields $\p_{x^i}$ and $\p_{\xi^\al}$ act on a function $f$ of the form
\eqref{f} by the rule
\begin{align}\p_{x^i}f=&\sum_{r=0}^m\sum_{\al_1<\cdots<\al_r}\p_{x^i}f_{\al_1...\al_r}\xi^{\al_1}\cdots\xi^{\al_r},\\
\p_{\xi^\al}f=&\sum_{r=1}^m\sum_{s=1}^r\sum_{\al_1<\cdots<\al_{s-1}<\al_s=\al<\al_{s+1}\cdots<\al_r}(-1)^{s-1}
f_{\al_1...\al_r}\xi^{\al_1}\cdots\xi^{\al_{s-1}}\xi^{\al_{s+1}}\cdots\xi^{\al_r}.\end{align} We will denote the vector
field $\p_{x^a}$ just by $\pa$.

Let $\M=(M,\O_\M)$ be a supermanifold and $\E$ a locally free sheaf of $\O_\M$-supermodules on $\M$, e.g. $\T_\M$. Let
$p|q$ be the rank of $\E$, then locally there exists a basis $(e_I,h_\Phi)_{I=1,...,p;\Phi=1,...,q}$ of sections of $\E$.
We denote such basis also by $e_A,$  where $e_{p+\Phi}=h_\Phi$. We will  always assume that $A,B,C=1,...,p+q$. For a point
$x\in M$ consider the vector space $\E_x=\E(V)/(\O_\M(V)_x\E(V))$, where $V\subset M$ is an open subset and $\O_\M(V)_x$
is the ideal in $\O_\M(V)$ consisting of functions vanishing at the point $x$. The vector space $\E_x$ does not depend on
choice of $V$; it is a real vector superspace of dimension $p|q$. For any open subset $V\subset M$ we have the projection
map from $\E(V)$ onto $\E_x$. For example, if $\E=\T_\M$, then $(\T_\M)_x$ is the tangent space $T_x\M$ to $\M$ at the
point $x$.

{\it A connection} on $\E$ is an even morphism $\na:\T_\M\otimes_\Real \E\to\E$ of sheaves of supermodules over $\Real$
such that $$\na_{fY}X=f\na_YX\quad\text{and}\quad\na_YfX=(Yf)X+(-1)^{|Y||f|}f\na_YX$$ for all homogeneous functions $f$,
vector fields $Y$ on $\M$ and sections $X$ of $\E$, here $|\cdot|\in\mathbb{Z}_2=\{{\bar 0},{\bar 1}\}$ denotes the
parity. In particular, $|\na_YX|=|Y|+|X|$. Locally we get the superfunctions $\Ga^A_{aB}$ such that
$\na_{\p_{a}}e_B=\Ga^A_{aB}e_A$. Obviously, $|\Ga^A_{aB}|=|a|+|A|+|B|$, where $|a|=|\p_{x^a}|$ and $|A|=|e_A|$.

The {\it curvature tensor} of the connection $\na$ is given by
\begin{equation}\label{R} R(Y,Z)=[\na_Y,\na_Z]-\na_{[Y,Z]},\end{equation}
where $Y$ and $Z$ are vector fields on $\M$. Let $\bar\na$ be a connection on $\T_\M$. Define the {\it covariant
derivatives} of $R$ with respect to $\bar\na$ as follows \begin{align}
\bar\na^r_{Y_r,...,Y_1}R(Y,Z)X=&\na_{Y_r}(\bar\na^{r-1}_{Y_{r-1},...,Y_1}R(Y,Z)X)-\bar\na^{r-1}_{\bar\na_{Y_r}Y_{r-1},...,Y_1}R(Y,Z)X\nonumber\\
&-(-1)^{|Y_r||Y_{r-1}|}\bar\na^{r-1}_{Y_{r-1},\bar\na_{Y_r}Y_{r-2},...,Y_1}R(Y,Z)X\nonumber
\\&-\cdots-(-1)^{|Y_r|(|Y_{r-1}|+\cdots+|Y_2|)}\bar\na^{r-1}_{Y_{r-1},...,Y_2,\bar\na_{Y_r}Y_{1}}R(Y,Z)X\label{naR}\\
&-(-1)^{|Y_r|(|Y_{r-1}|+\cdots+|Y_1|)}\bar\na^{r-1}_{Y_{r-1},...,Y_1}R(\bar\na_{Y_r}Y,Z)X\nonumber\\
&-(-1)^{|Y_r|(|Y_{r-1}|+\cdots+|Y_1|+|Y|)}\bar\na^{r-1}_{Y_{r-1},...,Y_1}R(Y,\bar\na_{Y_r}Z)X\nonumber\\
&-(-1)^{|Y_r|(|Y_{r-1}|+\cdots+|Y_1|+|Y|+|Z|)}\bar\na^{r-1}_{Y_{r-1},...,Y_1}R(Y,Z)\na_{Y_r}X,\nonumber
\end{align} here $r\geq 1$, $Y_r,...,Y_1,Y,Z\in \T_\M(M)$ are homogeneous and $X\in\E(M)$.
We assume that $\bar\na^0R=R$. It holds $|\bar\na^r_{Y_r,...,Y_1}R(Y,Z)|=|Y_r|+\cdots+|Y_1|+|Y|+|Z|$. If $\E=\T_\M$ and
$\bar\na=\na$, then we get the usual covariant derivatives of $R$.

%Consider local coordinates $(x^a)$ over $U\subset M$ and local basis $(e_A)$ of sections of $\E(U)$.
Let $r\geq 0$. Define the components $\bar\na^r_{a_r,...,a_1}R^A_{Bab}$ of $\bar\na^rR$ by the condition
$$\bar\na^r_{\p_{a_r},...,\p_{a_1}}R(\p_a,\p_b)e_B=\bar\na^r_{a_r,...,a_1}R^A_{Bab}e_A.$$ Then,
$|\bar\na^r_{a_r,...,a_1}R^A_{Bab}|=|a_r|+\cdots+|a_1|+|a|+|b|+|A|+|B|$. It is easy to show the following
\begin{equation}\label{Rcoord}
R^A_{Bab}=\pa\Ga^A_{bB}+(-1)^{|a|(|b|+|B|+|C|)}\Ga^C_{bB}\Ga^A_{aC}-(-1)^{|a||b|}(\pb\Ga^A_{aB}+(-1)^{|b|(|a|+|B|+|C|)}\Ga^C_{aB}\Ga^A_{bC})\end{equation}
and
\begin{align}\bar\na^r_{a_r,...,a_1}R^A_{Bab}=&\p_{a_r}(\bar\na^{r-1}_{a_{r-1},...,a_1}R^A_{Bab})
+(-1)^{|a_r|(|a_{r-1}|+\cdots+|a_1|+|a|+|b|+|B|+|C|)} \bar\na^{r-1}_{a_{r-1},...,a_1}R^C_{Bab}\Ga^A_{a_rC}\nonumber\\
&-\bar\Ga^c_{a_ra_{r-1}}\bar\na^{r-1}_{c,a_{r-2},...,a_1}R^A_{Bab}
-(-1)^{(|c|+|a_{r-2}|)|a_{r-1}|}\bar\Ga^c_{a_ra_{r-2}}\bar\na^{r-1}_{a_{r-1},c,a_{r-3},...,a_1}R^A_{Bab}\nonumber\\
&-\cdots-(-1)^{(|c|+|a_{1}|)(|a_{r-1}|+\cdots|a_2|)}\bar\Ga^c_{a_ra_{1}}\bar\na^{r-1}_{a_{r-1},...,a_2,c}R^A_{Bab}\label{naRcoord}\\
&-(-1)^{(|c|+|a|)(|a_{r-1}|+\cdots+|a_1|)}\bar\Ga^c_{a_ra}\bar\na^{r-1}_{a_{r-1},...,a_1}R^A_{Bcb}\nonumber\\
&-(-1)^{(|c|+|b|)(|a_{r-1}|+\cdots+|a_1|+|a|)}\bar\Ga^c_{a_rb}\bar\na^{r-1}_{a_{r-1},...,a_1}R^A_{Bac}\nonumber\\
&-(-1)^{(|C|+|B|)(|a_{r-1}|+\cdots+|a_1|+|a|+|b|)}\Ga^C_{a_rB}\bar\na^{r-1}_{a_{r-1},...,a_1}R^A_{Cab}.\nonumber
\end{align}

\section{Holonomy of smooth manifolds}\label{smoothhol}

In this section we recall some standard facts about holonomy of connections on vector bundles over smooth manifolds, see
e.g. \cite{K-N,Be}.

Let $E$ be a vector bundle over a connected smooth manifold $M$ and $\na$ a connection on $E$. It is known that for any
smooth curve $\gamma: [a,b]\subset\Real\to M$ and any  $X_0\in E_{\gamma(a)}$ there exists a unique section $X$ of $E$
defined along the curve $\gamma$ and satisfying the differential equation $\nabla_{\dot\gamma(s)}X=0$ with the initial
condition $X_{\gamma(a)}=X_0$. Consequently, for any  smooth curve $\gamma: [a,b]\subset\Real\to M$ we obtain the
isomorphism $\tau_{\gamma}:E_{\gamma(a)}\to E_{\gamma(b)}$ defined by $\tau_{\gamma}:X_0\mapsto X_{\gamma(b)}$. The
isomorphism $\tau_{\gamma}$ is called {\it the parallel displacement along the curve} $\gamma$. The parallel displacement
can be defined in the obvious way also for piecewise smooth curves.
 Let $x\in M$. {\it The holonomy group} $\Hol(\na)_x$ of the connection $\na$ at the point $x$ is the subgroup of $\GL(E_x)$
that consists  of parallel displacements  along all piecewise smooth loops at the point $x\in M$. If we consider only
null-homotopic loops, we get {\it the restricted holonomy group} $\Hol(\na)^0_x$. Obviously, $\Hol(\na)^0_x\subset
\Hol(\na)_x$ is a subgroup. If the manifold $M$ is simply connected, then $\Hol(\na)^0_x=\Hol(\na)_x$. It can be proved
that {\it the group $\Hol(\na)_x$ is a Lie subgroup of the Lie group $\GL(E_x)$ and the group $\Hol^0_x$ is the connected
identity component of the Lie group $\Hol(\na)_x$.} The Lie algebra $\hol(\na)_x$ of the Lie group $\Hol(\na)_x$ (and of
$\Hol(\na)^0_x$) is called {\it the holonomy algebra} of the connection $\na$ at the point $x$. Since the manifold $M$ is
connected,  the holonomy groups of $\na$ at different points of $M$ are isomorphic.

Remark that by the  holonomy group (resp. holonomy algebra) we understand not just the Lie group $\Hol(\na)_x$  (resp. Lie
algebra $\hol(\na)_x$), but the Lie group $\Hol(\na)_x$ with the representation $\Hol(\na)_x\hookrightarrow \GL(E_x)$
(resp. the Lie algebra $\hol(\na)_x$ with the representation $\hol\hookrightarrow\gl(E_x)$). These representations are
called {\it the holonomy representations}.

The theorem of Ambrose and Singer states that {\it the holonomy algebra $\hol(\na)_x$ coincides with the vector subspace
of $\gl(E_x)$ spanned by the elements of the form $$\tau_\ga^{-1}\circ R_y(Y,Z)\circ\tau_\ga,$$ where $R$ is the curvature
tensor of the connection $\na$, $\gamma$ is any curve in $M$ beginning at the point $x$; $y\in M$ is the end-point of the
curve $\gamma$ and $Y,Z\in T_yM$.}

Note that if $E=TM$ is the tangent bundle of $M$, then $$\tau_\ga^{-1}\circ
\na^r_{{Y_r},...,{Y_1}}R_y(Y,Z)\circ\tau_\ga\in\hol(\na)_x,$$ where $Y,Z,Y_1,...,Y_r\in T_yM$.

A section $X\in\Ga(E)$ is called {\it parallel} if $\na X=0$. This is equivalent to the condition that $X$ is parallel
along all curves in $M$, i.e. for any curve $\gamma:[a,b]\to M$ holds $\tau_\gamma X_{\ga(a)}=X_{\ga(b)}$.

The following  theorem is one of the main applications of the holonomy.

\begin{theorem}\label{paralsecprost}
Let $M$ be a smooth manifold, $E$ a vector bundle over $M$  and $\na$ a connection on $E$.
 Then the following conditions are equivalent:
\begin{description} \item[(i)] there exists a non-zero parallel section $X\in \Ga(E)$;
\item[(ii)] for any point $x\in M$ the holonomy group $\Hol(\na)_x$ preserves a non-zero vector $X_x\in E_x$.
%\item[(iii)] there exists a non-zero section $X\in \Ga(E)$ that is parallel along all curves in $M$.
\end{description}
%If the manifold  $M$ is simply connected, then  condition (ii) can be written as follows
%\begin{description}
%\item[(ii)] for any point $x\in M$ the holonomy algebra $\hol(\na)_x$ annihilates a non-zero vector $X_x\in E_x$.
%\end{description}
\end{theorem}

{\it To prove} the inclusion (i)$\Rightarrow$(ii) in Theorem \ref{paralsecprost} it is enough to take the value $X_x\in
E_x$. Since $X$ is invariant under the parallel displacements, the vector $X_x$ is invariant under the parallel
displacements along the loops at the point $x$, i.e.
 under  the holonomy representation.
Conversely, for a given vector $X_x\in E_x$ define the section $X\in\Ga(E)$ such that at any point  $y\in M$ holds
$X_y=\tau_{\gamma}X_x$,
 where $\gamma$ is any curve beginning at  $x$ and ending at $y$. From  condition (ii) it follows that $X_y$ does not depend
 on the choice of the curve $\gamma$. $\Box$

A vector subbundle  $F\subset E$ is called {\it parallel} if for all $Y\in\Ga(TM)$ and all $X\in\Ga(F)$ holds $\na_Y
X\in\Ga(F)$. This is equivalent to the condition that $F$ is parallel along all curves in $M$, i.e. for any curve
$\gamma:[a,b]\to M$ holds $\tau_\gamma F_{\ga(a)}=F_{\ga(b)}$.

\begin{theorem}\label{paralsubbund}
Let $M$ be a smooth manifold, $E$ a vector bundle over $M$  and $\na$ a connection on $E$.
 Then the following conditions are equivalent:
\begin{description} \item[(i)] there exists a  parallel vector subbundle  $F\subset E$ of rank $p$;
\item[(ii)] for any point $x\in M$ the holonomy group $\Hol(\na)_x$ preserves a  vector subspace $F_x\subset E_x$ of dimension
$p$.
\end{description}
\end{theorem}

\section{Definition of the holonomy for a connection on a sheaf over a supermanifold}\label{superhol}

Let $(\M,\O_\M)$ be a supermanifold, $\E$ a locally free sheaf of $\O_\M$-supermodules of rank $p|q$ on $\M$ and $\na$ a
connection on $\E$.
 Consider the  vector  bundle $E$ over $M$ defined as $E=\cup_{x\in M} \E_x$. The rank  of $E$ is $p+q$.
 Define the  subbundles $E_{\bar 0}=\cup_{x\in M} (\E_x)_{\bar 0}$ and $E_{\bar 1}=\cup_{x\in M} (\E_x)_{\bar 1}$ of $E$.
 Obviously, the restriction $\t\na=(\na|_{\Ga(TM)\otimes \Ga(E)})^\ti:\Ga(TM)\otimes \Ga(E)\to\Ga(E)$ is a connection on $E$.
Since $\na$ is even, the subbundles $E_{\bar 0},E_{\bar 1}\subset E$ are parallel.
 Let $\ga:[a,b]\subset\Real\to M$ be a curve and $\tau_\gamma:E_{\ga(a)}\to E_{\ga(b)}$ the parallel displacement
along $\ga$. Since the subbundles $E_{\bar 0},E_{\bar 1}\subset E$ are parallel, we have $\tau_\gamma (E_{\bar
0})_{\ga(a)}= (E_{\bar 0})_{\ga(b)}$ and $\tau_\gamma (E_{\bar 1})_{\ga(a)}= (E_{\bar 1})_{\ga(b)}$. We get the even
isomorphism $$\tau_\gamma:\E_{\ga(a)}\to \E_{\ga(b)}$$ of vector superspaces. We  call this isomorphism {\it the parallel
displacement} in $\E$ along $\ga$.

\begin{remark}  In \cite{Go} a parallel displacement in $\T_\M$ along  supercurves $\gamma:\Real^{1|1}\to\M$ is
introduced. The proof of the existence of the parallel displacement (p. 14) shows that the parallel displacements in
$\T_\M$ along a supercurve $\gamma:\Real^{1|1}\to\M$ and along the underlying curve $\t\gamma:\Real\to M$ coincide and
they coincide with our definition for the case $\E=\T_\M$.
%Parallel transport along supercurves has been carefully studied in \cite{Dumitr}.
\end{remark}

\begin{definit}\label{defhol}
Let $(\M,\O_\M)$ be a supermanifold, $\E$ a locally free sheaf of $\O_\M$-supermodules on $\M$ and $\na$ a connection on
$\E$.  The holonomy algebra $\hol(\na)_x$ of the connection $\na$ at a point $x\in M$ is the supersubalgebra of the Lie
superalgebra $\gl(\E_x)$ generated by the operators of the form $$\tau_\ga^{-1}\circ
\bar\na^r_{{Y_r},...,{Y_1}}R_y(Y,Z)\circ\tau_\ga:\E_x\to\E_x,$$ where $\gamma$ is any curve in $M$ beginning at the point
$x$; $y\in M$ is the end-point of the curve $\gamma$, $r\geq 0$, $Y,Z,Y_1,...,Y_r\in T_y\M$ and $\bar\na$ is a connection
on $\T_\M|_U$ for an open neighbourhood $U\subset M$ of $y$.
\end{definit}

\begin{prop}\label{correct} The definition of the holonomy algebra $\hol(\na)_x$ does not depend on the choice of the
connection $\bar\na$.\end{prop}

{\bf Proof.}  Let $\gamma$ be a curve in $M$ beginning at the point $x$ and ending at a point $y\in M$. Let $U\subset M$
be an open neighbourhood of the point $y$. For any connection $\hat\na$ on $\T_\M|_{U}$ and any integer $t\geq 0$ define
the vector space $$L(\hat\na)_t=\spa\{\hat\na^r_{{Y_r},...,{Y_1}}R_y(Y,Z)|0\leq r\leq t,\,\,Y,Z,Y_1,...,Y_r\in T_y\M\}.$$
Clearly, $L(\hat\na)_t$ does not depend on the choice of $U$.   Let $(U,x^a)$ be a system of local coordinates such that
$y\in U$. Let $\bar\na$ be a connection on $\T_\M|_{U}$.
 Denote by $\vec\na$  the connection on $\T_\M|_{U}$ such that $\vec\na\pa=0$. To prove the proposition  it is  enough to
 show that  for any $t\geq 0$ we  have $L(\bar\na)_t=L({\vec\na})_t$. This will follow from the following lemma.

\begin{lem}\label{lemcorrect} For any $t\geq 0$ it holds
$$\bar\na^t_{\p_{a_t},...,\p_{a_1}}R(\p_{a_0},\p_{a_{-1}})=\vec\na^t_{\p_{a_t},...,\p_{a_1}}R(\p_{a_0},\p_{a_{-1}})
+\sum_{s=0}^{t-1}\sum_{(b_s,...,b_{-1})}
B_{a_t...a_{-1},b_s...b_{-1}}\vec\na^s_{\p_{b_s},...,\p_{b_1}}R(\p_{b_0},\p_{b_{-1}}),$$ where
$B_{a_t...a_{-1},b_s...b_{-1}}\in\O_\M(U)$.
\end{lem}
{\it Proof.} We will prove the lemma by the induction over $t$. For $t=0$ there is nothing to prove. Fix $t>0$. Suppose
that the lemma is true for all $r<t$ and prove it for $r=t$.

Using \eqref{naRcoord} and the induction hypothesis, we get\\ $\bar\na^r_{\p_{a_r},...,\p_{a_1}}R(\p_{a_0},\p_{a_{-1}})$
\begin{align*} =&[\na_{\p_r},\bar\na^{r-1}_{\p_{a_{r-1}},...,\p_{a_1}}R(\p_{a_0},\p_{a_{-1}})]\\
&-\sum_{l=-1}^{r-1}(-1)^{(|c|+|a_l|)(|a_{r-1}|+\cdots+|a_{l+1}|)}
\bar\Ga^c_{a_ra_{l}}\bar\na^{r-1}_{\p_{a_{r-1}},...,\p_{a_{l+1}},\p_c,\p_{a_{l-1}},...,\p_{a_1}}R(\p_{a_0},\p_{a_{-1}})\\
=&[\na_{\p_r},\vec\na^{r-1}_{\p_{a_{r-1}},...,\p_{a_1}}R(\p_{a_0},\p_{a_{-1}})+\sum_{s=0}^{r-2}\sum_{(b_s,...,b_{-1})}
B_{a_{r-1}...a_{-1},b_s...b_{-1}}\vec\na^{s}_{\p_{b_s},...,\p_{b_1}}R(\p_{b_0},\p_{b_{-1}})]\\
&-\sum_{l=-1}^{r-1}(-1)^{(|c|+|a_l|)(|a_{r-1}|+\cdots+|a_{l+1}|)}
\bar\Ga^c_{a_ra_{l}}\bar\na^{r-1}_{\p_{a_{r-1}},...,\p_{a_{l+1}},\p_c,\p_{a_{l-1}},...,\p_{a_1}}R(\p_{a_0},\p_{a_{-1}})\\
%\end{align*}\begin{align*}
=&[\na_{\p_r},\vec\na^{r-1}_{\p_{a_{r-1}},...,\p_{a_1}}R(\p_{a_0},\p_{a_{-1}})]+
\sum_{s=0}^{r-2}\sum_{(b_s,...,b_{-1})}\p_{a_r}(B_{a_{r-1}...a_{-1},b_s...b_{-1}})\vec\na^{s}_{\p_{b_s},...,\p_{b_1}}
R(\p_{b_0},\p_{b_{-1}})\\ &+\sum_{s=0}^{r-2}\sum_{(b_s,...,b_{-1})}(-1)^{|a_r||B_{a_{r-1}...a_{-1},b_s...b_{-1}}|}
B_{a_{r-1}...a_{-1},b_s...b_{-1}}[\na_{\p_{a_r}},\vec\na^{s}_{\p_{b_s},...,\p_{b_1}}R(\p_{b_0},\p_{b_{-1}})]\\
&-\sum_{l=-1}^{r-1}(-1)^{(|c|+|a_l|)(|a_{r-1}|+\cdots+|a_{l+1}|)}
\bar\Ga^c_{a_ra_{l}}\bar\na^{r-1}_{\p_{a_{r-1}},...,\p_{a_{l+1}},\p_c,\p_{a_{l-1}},...,\p_{a_1}}R(\p_{a_0},\p_{a_{-1}})\\
%\end{align*}\begin{align*}
=&\vec\na^{r}_{\p_{a_{r}},...,\p_{a_1}}R(\p_{a_0},\p_{a_{-1}})+
\sum_{s=0}^{r-2}\sum_{(b_s,...,b_{-1})}\p_{a_r}(B_{a_{r-1}...a_{-1},b_s...b_{-1}})\vec\na^{s}_{\p_{b_s},...,\p_{b_1}}
R(\p_{b_0},\p_{b_{-1}})\\ &+\sum_{s=0}^{r-2}\sum_{(b_s,...,b_{-1})}(-1)^{|a_r||B_{a_{r-1}...a_{-1},b_s...b_{-1}}|}
B_{a_{r-1}...a_{-1},b_s...b_{-1}}\vec\na^{s+1}_{\p_{a_r}\p_{b_s},...,\p_{b_1}}R(\p_{b_0},\p_{b_{-1}})\\
&-\sum_{l=-1}^{r-1}(-1)^{(|c|+|a_l|)(|a_{r-1}|+\cdots+|a_{l+1}|)}
\bar\Ga^c_{a_ra_{l}}\bar\na^{r-1}_{\p_{a_{r-1}},...,\p_{a_{l+1}},\p_c,\p_{a_{l-1}},...,\p_{a_1}}R(\p_{a_0},\p_{a_{-1}}).
\end{align*} The proof of the lemma follows from the induction hypothesis applied to the last term. $\Box$

 The proposition is proved. $\Box$

The next proposition simplifies the expression for the holonomy algebra. In particular it shows that it is not necessary
to take the covariant derivatives of the curvature tensor in the directions of the vectors tangent to $M$.

\begin{prop} The holonomy algebra $\hol(\na)_x$ coincides with the supersubalgebra of the Lie
superalgebra $\gl(\E_x)$ generated by the operators of the form $$\tau_\ga^{-1}\circ
R_y(\p_i,\p_j)\circ\tau_\ga\qquad\text{and}\qquad\tau_\ga^{-1}\circ
\bar\na^r_{\p_{\ga_r},...,\p_{\ga_1}}R_y(\p_\ga,\p_a)\circ\tau_\ga,$$ where $\gamma$ is any curve starting at the point
$x$; $y$ is the end-point of the curve $\ga$, $r\geq 0$, $\ga_r>\cdots>\ga_1$, $x^a$ are local coordinates on $\M$ over an
open neighbourhood $U$ of the point $y$ and $\bar\na$ is a connection on $\T_\M|_U$.
\end{prop}

{\bf Proof.} Let $\g$ be the supersubalgebra of $\gl(\E_x)$ generated by the operators as in the formulation of the
proposition without the assumption $\ga_r>\cdots>\ga_1$. We will prove that $\hol(\na)_x=\g$. Fix a curve $\gamma$
starting at the point $x$. Let $y$ be the end-point of $\gamma$, $U\subset M$ an open neighbourhood of $y$, $x^a$ local
coordinates on $\M$ over $U$ and $\bar\na$ a connection on $\T_\M|_U$. It is enough to show that any element
$\eta=\tau_\ga^{-1}\circ \bar\na^r_{\p_{a_r},...,\p_{a_1}}R_y(\p_a,\p_b)\circ\tau_\ga\in\hol(\na)_x$ belongs to $\g$. We
will prove this statement by the induction over $r$. By Lemma \ref{lemcorrect}, we may assume that the connection
$\bar\na$ is flat such that $\na\p_a=0$. For $r=0$ there is nothing to prove. Fix $t>0$ suppose that the statement is true
for all $r<t$ and prove it for $r=t$.

\begin{lem} Let $r\geq 2$ and $Y_1,...,Y_r,Y,Z\in \T_\M(U)$ be matually commuting homogeneous parallel vector fields. Then
\begin{multline}\label{formzam}
\bar\na^r_{Y_r,...,Y_1}R(Y,Z)=(-1)^{|Y_s||Y_{s-1}|}\bar\na^r_{Y_r,...Y_{s+1},Y_{s-1},Y_s,Y_{s-2},...,Y_1}R(Y,Z)
\\+ \sum_{I\subset \{s+2,...,r\}} (-1)^{S(I)}\Big[\bar\na^{[I]}_{Y_{i_1},...,Y_{i_{[I]}}}R(Y_{s+1},Y_s),
\bar\na^{r-2-[I]}_{Y_{i'_1},...,Y_{i'_{[I']}},Y_{s-1},...,Y_1}R(Y,Z)\Big],
\end{multline}
 where
$I=\{i_1,...,i_{[I]}\}\subset \{s+2,...,r\}$ is a subset with $i_1>\cdots >i_{[I]}$, $[I]$ is the number of elements in
$I$, $I'=\{i'_1,...,i'_{[I']}\}=\{s+2,...,r\}\backslash I$, $i'_1>\cdots> i'_{[I']}$ and $S(I)$ is the sign of the
permutation $(i_1,...,i_{[I]},s+1,s,i'_1,...,i'_{[I']})$ of the numbers $r,...,1$.
\end{lem}

{\it Proof.} Using \eqref{naR}, the assumptions of the lemma  and the Jacobi super identity, we get
\begin{align}\bar\na^r_{Y_r,...,Y_1}R(Y,Z)=&[\na_{Y_r},[\na_{Y_{r-1}},...,[\na_{Y_{s+1}},[\na_{Y_s},\bar\na^{s-1}_{Y_{s-1},...,Y_1}R(Y,Z)]]...]]\nonumber
\\=&[\na_{Y_r},[\na_{Y_{r-1}},...,[[\na_{Y_{s+1}},\na_{Y_s}],\bar\na^{s-1}_{Y_{s-1},...,Y_1}R(Y,Z)]...]]\nonumber\\&+(-1)^{|Y_s||Y_{s-1}|}
[\na_{Y_r},[\na_{Y_{r-1}},...,[\na_{Y_{s}},[\na_{Y_{s+1}},\bar\na^{s-1}_{Y_{s-1},...,Y_1}R(Y,Z)]]...]]\nonumber\\=&
[\na_{Y_r},[\na_{Y_{r-1}},...,[R({Y_{s+1}},{Y_s}),\bar\na^{s-1}_{Y_{s-1},...,Y_1}R(Y,Z)]...]]\nonumber\\&+
(-1)^{|Y_s||Y_{s-1}|}\bar\na^r_{Y_r,...Y_{s+1},Y_{s-1},Y_s,Y_{s-2},...,Y_1}R(Y,Z).\nonumber
\end{align}
Applying the Jacobi super identity $r-s$ times to the first term of the last equality, we get the proof of the lemma.
$\Box$

Note that  for a curve $\mu(t)$ in $M$ such that $\mu(0)=y$ it holds
\begin{equation}\label{namu} \bar\na^r_{\dot\mu(0),\p_{a_{r-1}},...,\p_{a_1}}R_y(\p_a,\p_b)=
\frac{d}{dt}\Big|_{t=0}((\tau_\mu)^t_0)^{-1}\circ\bar\na^{r-1}_{\p_{a_{r-1}},...,\p_{a_1}}R_{\mu(t)}(\p_a,\p_b)\circ(\tau_\mu)^t_0,
\end{equation} where $(\tau_\mu)^t_0$ is the parallel displacement along the curve $\mu|_{[0,t]}$.
Fix an element $\eta=\tau_\ga^{-1}\circ \bar\na^r_{\p_{a_r},...,\p_{a_1}}R_y(\p_a,\p_b)\circ\tau_\ga\in\hol(\na)_x$.  If
$a_r\leq n$, then there exists a curve $\mu(t)$ with $\mu(0)=y$ and $\dot\mu(t)=\p_{a_r}(\mu(t))$. From \eqref{namu} and
the induction hypothesis it follows that $\eta\in\g$. Suppose that there exists $s$ such that $r-1\geq s\geq 1$ and
$a_s\leq n$. Let $s$ be maximal with this property. Applying to $\eta$ several times \eqref{formzam} and using the
induction hypothesis, we get $\eta-(-1)^{r-s}\tau_\ga^{-1}\circ
\bar\na^r_{\p_{a_s},\p_{a_{r}},...,\p_{a_{s+1}},\p_{a_{s-1}},...,\p_{a_{1}}}R_y(\p_a,\p_b)\circ\tau_\ga\in\g$. By the
above argumentation, $\eta\in\g$. Now we may assume that  $a_r,...,a_1>n$. If  $a\leq n$ or $b\leq n$, then using the
second Bianchi super identity we get $\eta=-(-1)^{|a_1||a|}\tau_\ga^{-1}\circ
\bar\na^r_{\p_{a_r},...,\p_{a_2},\p_a}R_y(\p_b,\p_{a_1})\circ\tau_\ga-\tau_\ga^{-1}\circ
\bar\na^r_{\p_{a_r},...,\p_{a_2},\p_b}R_y(\p_{a_1},\p_a)\circ\tau_\ga\in\g.$ Thus, $\hol(\na)_x=\g$. Equality
\eqref{formzam} shows that $\g$ coincides with  the supersubalgebra of $\gl(\E_x)$ generated by the operators as in the
formulation of the proposition. $\Box$

Let $E$ be the vector bundle over $M$ and $\t\na$ the connection on $E$ as above. Then the holonomy algebra
$\hol(\t\na)_x$ is contained in $(\hol(\na)_x)_{\bar 0}$, but these Lie algebras must not coincide, this shows  the
following example.

\begin{ex} Consider the supermanifold $\Real^{0|1}=(\{0\},\Lambda_\xi)$.
Define the connection $\na$ on $\T_{\Real^{0|1}}$ by $\na_{\p_{\xi}}\p_{\xi}=\xi\p_{\xi}.$ Then, $\hol(\t\na)_0=\{0\}$ and
$(\hol(\na)_0)_{\bar 0}=\hol(\na)_0=\gl(0|1)$.\end{ex}

Now we define the holonomy group.  Recall that a {\it Lie supergroup} $\mathcal G=(G,\O_\mathcal G)$ is a group object in
the category of supermanifolds. The underlying smooth manifold $G$ is a Lie group. The Lie superalgebra $\g$ of $\mathcal
G$ can be identified with the tangent space to $\mathcal G$ at the identity $e\in G$. The Lie algebra of the Lie group $G$
is the even part $\g_{\bar 0}$ of the Lie superalgebra $\g$.

Any Lie supergroup $\mathcal G$ is uniquely given by a pair $(G,\g)$ ({\it Harish-Chandra pair}), where $G$ is a Lie
group, $\g=\g_{\bar 0}\oplus\g_{\bar 1}$ is a Lie superalgebra such that $\g_{\bar 0}$ is the Lie algebra of the Lie group
$G$ and there exists a representation $\Ad$ of $G$ on $\g$ that extends the adjoint representation of $G$ on $\g_{\bar 0}$
and the differential of $\Ad$ coincides with the Lie superbracket of $\g$ restricted to $\g_{\bar 0}\times\g_{\bar 1}$,
see \cite{DelMor,Go}.

Denote by $\Hol(\na)^0_x$ the connected Lie subgroup of the Lie group $\GL((\E_x)_{\bar 0})\times \GL((\E_x)_{\bar 1})$
corresponding to the Lie subalgebra $(\hol(\na)_x)_{\bar 0}\subset \gl((\E_x)_{\bar 0})\oplus \gl((\E_x)_{\bar
1})\subset\gl(\E_x)$. Let $\Hol(\na)_x$ be the Lie subgroup of the Lie group $\GL((\E_x)_{\bar 0})\times \GL((\E_x)_{\bar
1})$ generated by the Lie groups $\Hol(\na)^0_x$ and $\Hol(\t\na)_x$. Clearly, the Lie algebra of the Lie group
$\Hol(\na)_x$ is $(\hol(\na)_x)_{\bar 0}$. Let $\Ad'$ be the representation of the connected Lie group $\Hol(\na)^0_x$ on
$\hol(\na)_x$ such that the differential of $\Ad'$ coincides with the Lie superbracket of $\hol(\na)_x$ restricted to
$(\hol(\na)_x)_{\bar 0}\times(\hol(\na)_x)_{\bar 1}$. Define the representation $\Ad''$ of the Lie group $\Hol(\t\na)_x$
on $\hol(\na)_x$ by the  rule $$\Ad''_{\tau_\mu}(\tau_\ga^{-1}\circ
\bar\na^r_{{Y_r},...,{Y_1}}R_y(Y,Z)\circ\tau_\ga)=\tau_{\mu}\circ \tau_\ga^{-1}\circ \bar\na^r_{ Y_r,..., Y_1}R_y( Y,
Z)\circ\tau_\ga\circ \tau_\mu^{-1}.$$ Note that $\Ad'|_{\Hol(\na)^0_x\cap \Hol(\t\na)_x}=\Ad''|_{\Hol(\na)^0_x\cap
\Hol(\t\na)_x}$. Consequently, we get a representation $\Ad$ of the group $\Hol(\na)_x$ on $\hol(\na)_x$. It is obvious
that $\Hol(\na)^0_x\cap \Hol(\t\na)_x=\Hol(\t\na)^0_x$ and if $M$ is simply connected, then $\Hol(\t\na)_x\subset
\Hol(\na)^0_x$ and $\Hol(\na)_x=\Hol(\na)^0_x$.

\begin{definit} The Lie supergroup $\Hols(\na)_x$ given by the Harish-Chandra pair $(\Hol(\na)_x,\hol(\na)_x)$ is called the
holonomy group of the connection $\na$ at the point $x$. The Lie supergroup $\Hols(\na)^0_x$ given by the Harish-Chandra
pair $(\Hol(\na)^0_x,\hol(\na)_x)$ is called the restricted holonomy group of the connection $\na$ at the point
$x$.\end{definit}

\section{Infinitesimal holonomy algebras}\label{sectinf}

In this section we define the infinitesimal holonomy algebra and  show that in the analytic case it coincides with the
holonomy algebra.

\begin{definit}
Let $(\M,\O_\M)$ be a supermanifold, $\E$ a locally free sheaf of $\O_\M$-supermodules on $\M$ and $\na$ a connection on
$\E$.  The infinitesimal holonomy algebra $\hol(\na)^{inf}_x$ of the connection $\na$ at a point $x\in M$ is the
supersubalgebra of the Lie superalgebra $\gl(\E_x)$ generated by the operators of the form
$$\bar\na^r_{{Y_r},...,{Y_1}}R_x(Y,Z),$$ where  $r\geq 0$, $Y,Z,Y_1,...,Y_r\in T_x\M$ and $\bar\na$ is a connection on
$\T_\M|_U$ for an open neighbourhood $U\subset M$ of $x$.
\end{definit}

Lemma \ref{lemcorrect} shows that the definition of  $\hol(\na)^{inf}_x$ does not depend on the choice of the connection
$\bar\na$.

\begin{theorem}
Let $\M=(M,\O_\M)$ be an analytic supermanifold, $\E$ a locally free sheaf of $\O_\M$-supermodules on $\M$ and $\na$ a
connection on $\E$.   Then $$\hol(\na)_x=\hol(\na)^{inf}_x.$$
\end{theorem}

{\bf Proof.} First we prove the following lemma.

\begin{lem} If $\gamma:[0,1]\subset\Real\to M$ is a piecewise analytic curve starting at $x\in M$ and ending at $y\in M$,
then $$\tau_\ga^{-1}\circ\hol(\na)^{inf}_y\circ\tau_\ga=\hol(\na)^{inf}_x.$$
\end{lem}

{\it Proof.} Consider some partition $0=t_0<t_1<\cdots<t_k=1$. Obviously, if the statement of the lemma holds for each
curve $\ga|_{[t_i,t_{i+1}]}$, then the lemma is true. Consequently we may assume that the image of $\ga$ is contained in
some coordinate neighbourhood $U\subset M$. Fix a connection $\bar\na$  on $\T_\M|_U$.

 Fix an element $\bar\na^r_{{\bar Y_r},...,{\bar Y_1}}R_y(\bar Y,\bar
Z)\in\hol(\na)^{inf}_y$. It is enough to prove that $\tau_\ga^{-1}\circ \bar\na^r_{{\bar Y_r},...,{\bar Y_1}}R_y(\bar
Y,\bar Z)\circ\tau_\ga\in\hol(\na)^{inf}_x$.  Let $Y, Z, Y_1,..., Y_r\in T_x\M$ be the vectors such that $\bar Y=\tau_\ga
Y,\bar Z=\tau_\ga Z,\bar Y_1=\tau_\ga Y_1,...,\bar Y_r=\tau_\ga Y_r$. For any $s,t\in [0,1]$ with $s\leq t$ denote by
$(\tau_\ga)^t_s:\E_{\ga(s)}\to\E_{\ga(t)}$ the parallel displacement along the curve $\ga|_{[s,t]}$. Consider the
endomorphism $$F(t)=((\tau_\ga)^t_0)^{-1}\circ \bar\na^r_{{(\tau_\ga)^t_0 Y_r},...,{(\tau_\ga)^t_0
Y_1}}R_{\ga(t)}((\tau_\ga)^t_0 Y,(\tau_\ga)^t_0 Z)\circ(\tau_\ga)^t_0:\E_x\to\E_x.$$ We must prove that
$F(1)\in\hol(\na)^{inf}_x$. Fix a basis of $\E_x$. Without loss of generality we my assume that the elements $F^A_B(t)$ of
the matrix of $F(t)$ are analytic functions of $t$, i.e. $F^A_B(t)=\sum_{k=0}^{\infty}F^A_{Bk}t^k$ for some real numbers
$F^A_{Bk}$. For each $k\geq 0$ denote by $F_k$ the endomorphism of $\E_x$ with the matrix $F^A_{Bk}$. Since
$F^A_{B0}=F^A_B(0)$, we have $F_0=F(0)=\bar\na^r_{{Y_r},...,{ Y_1}}R_{x}( Y, Z)\in\hol(\na)^{inf}_x$. Further,
$$\frac{d}{dt}F(t)=\lim_{s\to 0}\frac{F(t+s)-F(t)}{s}=((\tau_\ga)^t_0)^{-1}\circ \bar\na^{r+1}_{{(\tau_\ga)^t_0
Y_{r+1}},{(\tau_\ga)^t_0 Y_r},...,{(\tau_\ga)^t_0 Y_1}}R_{\ga(t)}((\tau_\ga)^t_0 Y,(\tau_\ga)^t_0 Z)\circ(\tau_\ga)^t_0,$$
where $Y_{r+1}=((\tau_\ga)^t_0)^{-1}\dot\ga(t)$. Consequently, $F_1=\frac{d}{dt}F(0)\in\hol(\na)^{inf}_x$. Similarly for
each $k>1$ $F_k=\frac{d^k}{(dt)^k}F(0)\in\hol(\na)^{inf}_x$. Thus, $F(1)=\sum_{k=0}^{\infty}F_{k}\in\hol(\na)^{inf}_x$.
The lemma is proved. $\Box$

Let now $\tau_\ga^{-1}\circ\bar \na^r_{{ Y_r},...,{ Y_1}}R_y( Y, Z)\circ\tau_\ga\in\hol(\na)_x$, where $\ga$ is any
piecewise smooth curve beginning at $x\in M$ and ending at $y\in M$. Let $\mu$ be a piecewise analytic curve beginning at
$x\in M$ and ending at $y\in M$ such that the loop $\ga*\mu^{-1}$ is null-homotopic. We have
\begin{multline*}\tau_\ga^{-1}\circ\bar
\na^r_{{ Y_r},...,{ Y_1}}R_y( Y, Z)\circ\tau_\ga=\tau_\mu^{-1}\circ\Big((\tau_\mu\circ\tau_\ga^{-1})\circ\bar \na^r_{{
Y_r},...,{ Y_1}}R_y( Y, Z)\circ(\tau_\ga\circ\tau_\mu^{-1})\Big)\circ\tau_\mu\\
=\tau_\mu^{-1}\circ\Big(\Ad''_{\tau_{\mu*\ga^{-1}}}\bar \na^r_{{ Y_r},...,{ Y_1}}R_y(Y,Z)\Big)\circ\tau_\mu.
\end{multline*}

Here $\tau_{\mu*\ga^{-1}}\in \Hol(\tilde\na)^0_y$ acts on the element $\eta=\bar \na^r_{{ Y_r},...,{ Y_1}}R_y( Y,
Z)\in\hol(\na)^{inf}_y\subset\hol(\na)_y$ as it was defined above. Note that the connection $\tilde\na$ is analytic. By
the classical result, the group $\Hol(\tilde\na)^0_y$ coincides with the infinitesimal holonomy group
$\Hol(\tilde\na)^{inf}_y$. Consequently, $\tau_{\mu*\ga^{-1}}=\exp\xi$ for some
$\xi\in\hol(\tilde\na)^{inf}_y\subset\hol(\na)^{inf}_y$. Finally, $$\Ad''_{\exp\xi}\eta=(\exp (d\Ad'')_\xi)\eta=(\exp
\ad_\xi)\eta=\sum_{k=0}^{\infty}\frac{\ad_\xi^k}{k!}\eta\in\hol(\na)_y^{inf}.$$ By the above lemma,
$\tau_\ga^{-1}\circ\bar \na^r_{{ Y_r},...,{ Y_1}}R_y( Y, Z)\circ\tau_\ga\in\hol(\na)_x^{inf}.$ Thus,
$\hol(\na)_x\subset\hol(\na)_x^{inf}$. The inverse inclusion is trivial. The theorem is proved. $\Box$

\section{Parallel sections}\label{parallel}
Let $(\M,\O_\M)$ be a supermanifold, $\E$ a locally free sheaf of $\O_\M$-supermodules on $\M$ and $\na$ a connection on
$\E$. A section $X\in\E(M)$ is called {\it parallel} if  $\na X=0$. Let $E$ be the vector bundle on $M$ and $\t\na$  the
connection on $E$ as above. For the section $\t X\in\Ga(M)$ we get $\t\na\t X=0.$ Hence for any curve
$\gamma:[a,b]\subset\Real\to M$ we have $\tau_\ga \t X_{\ga(a)}=\t X_{\ga(b)}$. Consequently, $\tau_\ga
X_{\ga(a)}=X_{\ga(b)}$, where $\tau_\ga:\E_{\ga(a)}\to\E_{\ga(b)}$, i.e. $X$ is {\it parallel along curves in $M$}.

Consider a system of local coordinates $(U,x^a)$ and a basis $e_A$ of $\E(U)$. Let $X\in\E(U)$, then $$X=X^Ae_A,\quad
\na_\pa X=\pa X^Ae_A+(-1)^{|a||X^A|}X^A\Ga^B_{aA}e_B.\footnote{We assume that $(-1)^{|a||X^A|}X^A\Ga^B_{a
A}=(-1)^{|a|}X_{\bar 0}^A\Ga^B_{a A}-(-1)^{|a|}X_{\bar1}^A\Ga^B_{\ga A}$, where $X^A=X^A_{\bar0}+X^A_{\bar1}$ is the
decomposition of $X^A$ in the sum of the even and odd parts.}$$ Thus the condition $\na X=0$ is equivalent to the
following condition in local coordinates
\begin{equation}\label{A+B} \pa X^A+(-1)^{|a||X^B|}X^B\Ga^A_{a B}=0\end{equation}
or to the conditions
\begin{align} \pai X^A+X^B\Ga^A_{iB}&=0,\label{A}\\
 \pga X^A+(-1)^{|X^B|}X^B\Ga^A_{\ga B}&=0.\label{B}\end{align}
 Equations \eqref{A} and \eqref{B} are equivalent to
\begin{align} (\p_{\ga_r}...\p_{\ga_1}(\pai X^A+X^B\Ga^A_{iB}))^\ti&=0,\label{A1}\\
 (\p_{\ga_r}...\p_{\ga_1}(\pga X^A+(-1)^{|X^B|}X^B\Ga^A_{\ga B}))^\ti&=0,\label{B1}\end{align}
where $0\leq r\leq m$. These equations can be written as
\begin{align} \pai
X_{\ga_1...\ga_r}^A+\sign_{\ga_1,...,\ga_r}\sum_{l=0}^r\sum_{\begin{smallmatrix}\{\al_1,...,\al_r\}=\{\ga_1,...,\ga_r\}\\\al_1<\cdots<\al_l,\al_{l+1}<\cdots<\al_r\end{smallmatrix}}\sign_{\al_1,...,\al_r}
X^B_{\al_1...\al_l}\Ga^A_{iB\al_{l+1}...\al_r}&=0,\label{A2}\\
 X_{\ga\ga_1...\ga_r}^A+\sign_{\ga_1,...,\ga_r}\sum_{l=0}^r\sum_{\begin{smallmatrix}\{\al_1,...,\al_r\}=
 \{\ga_1,...,\ga_r\}\\
 \al_1<\cdots<\al_l,\al_{l+1}<\cdots<\al_r\end{smallmatrix}}\sign_{\al_1,...,\al_r}(-1)^l
X^B_{\al_1...\al_l}\Ga^A_{\ga B\al_{l+1}...\al_r}&=0.\label{B2}\end{align}

Using this, we can prove the following proposition.

\begin{prop}\label{prop1} Let $\M=(M,\O_\M)$ be a supermanifold, $\E$ a locally free sheaf of $\O_\M$-supermodules on $\M$ and $\na$ a
 connection on $\E$.
Then   a parallel section $X\in\E(M) $ is uniquely defined by its value at any point $x\in M$.\end{prop} {\it Proof.} Let
$\na X=0$, $x\in M$ and $X_x$ be the value of $X$ at the point $x$. Since $X$ is parallel along curves in $M$, using
$X_x$, we can find the values of $X$ at all points of $M$. Consider the local coordinates as above. As we know the values
of $X$ at all points, we know the functions $\t X^A$. Using \eqref{B2} for $r=0$, we can find the functions $X^A_\ga$.
Namely, $X^A_\ga=-\t X^B\t\Ga_{\ga B}^A$. Using \eqref{B2} for $r=1$, we get $X^A_{\ga\ga_1}=-\t X^B\Ga^A_{\ga
B\ga_1}+X^B_{\ga_1}\t\Ga^A_{\ga B}$. In the same way we can find all functions $X^A_{\ga\ga_1...\ga_r}$, i.e. we know the
functions $X^A$ and we reconstruct the section $X$ in any coordinate system. The proposition is proved. $\Box$

\begin{theorem}\label{paralsec}
Let $\M=(M,\O_\M)$ be a supermanifold, $\E$ a locally free sheaf of $\O_\M$-supermodules on $\M$ and $\na$ a connection on
$\E$.   Then the following conditions are equivalent:
\begin{description} \item[(i)] there exists a non-zero paralel section $X\in\E(M)$;

\item[(ii)] for any point $x\in M$ there exists a non-zero vector $X_x\in \E_x$ annihilated by the holonomy algebra $\hol(\na)_x$
and preserved by the group $\Hol(\t\na)_x$.
\end{description}
%If the manifold  $M$ is simply connected, then  condition (ii) can be written as follows
%\begin{description}
%\item[(ii)] for any point $x\in M$ there exists a non-zero vector $X_x\in \E_x$ annihilates by the holonomy algebra $\hol(\na)_x$.
%\end{description}
\end{theorem}

{\bf Proof of Theorem \ref{paralsec}.}  Suppose that a section $X\in\E(M)$ is parallel, i.e. $\na X=0$. Let $U\subset M$
be an open subset and  $\bar\na$ a connection on $\T_\M|_U$. From \eqref{naR} it follows that
\begin{equation}\label{eq1}\bar\na^r_{{Y_r},...,{Y_1}}R(Y,Z)X=0\end{equation} for any vector fields $Y,Z,Y_1,...,Y_r\in\T_\M(U)$.
Since $X$ is parallel, for any curve $\ga$ in $M$ beginning at the point $x\in M$ and ending at a point $y\in M$, we
have $X_y=\tau_\ga X_x$. Hence to prove the implication (i)$\Rightarrow$(ii) it is enough to consider  \eqref{eq1} at the
point $y$.

Let us prove the implication (ii)$\Rightarrow$(i). Let $x\in M$ be any point and suppose there exists a non-zero vector
$X_x\in \E_x$ annihilated by the holonomy algebra $\hol(\na)_x$ and preserved by the group $\Hol(\t\na)_x$. From Theorem
\ref{paralsecprost} it follows that there exists a section $X_0\in\Ga(E)$ such that $\t\na X_0=0$ and $(X_0)_x=X_x$. Fix a
coordinate neighborhood $(U,x^a)$ on $\M$ and a local basis $e_A$ of $\E(U)$. Then $\tilde e_A$ is a local basis of
$\Ga(U,E)$ and we get the functions $X_0^A\in \O_M(U)$ such that $X_0=X_0^A\tilde e_A$ on $U$. Using \eqref{B2} and
$X^A_0$, as in the proof of Proposition \ref{prop1}, define functions $X^A_{\ga\ga_1...\ga_r}\in\O_M(U)$ for all
$\ga<\ga_1<\cdots<\ga_r$, $0\leq r\leq m-1$. This gives us functions $X^A\in\O_\M(U)$ such that $\t X^A=X_0^A$. Consider
the section $X=X^Ae_A\in\E(U)$. We claim that $\na X=0$. To prove this it is enough to show that the functions $X^A$
satisfy \eqref{A1} and \eqref{B1} for all $\ga_1<\cdots<\ga_r$, $0\leq r\leq m$ and any $\ga$, then $X^A$ will satisfy
\eqref{A1} and \eqref{B1} for all $\ga_1,...,\ga_r$ and $\ga$. Note that, by construction, the functions $X^A$ satisfy
\eqref{B1} for $\ga<\ga_1<\cdots<\ga_r$.

For the proof we use the induction over $r$. Parallel to this we will prove that
\begin{align}\label{st2} (\p_{a_r}...\p_{a_{s+1}}((-1)^{(|A|+|B|)|X^B|}\bar\na^{s-1}_{a_s,...,a_2}
R^A_{Ba_1a}X^B))^\ti&=0\quad\text{for all } r\geq 1 \text{ and }  1\leq s\leq r. \end{align}

For $r=0$ \eqref{A1} follows from the fact that $\t\na X_0=0$; \eqref{B1} follows from the definition of functions
$X^A_\ga$. Fix $r_0>0$. Suppose that \eqref{A1} and \eqref{B1} hold for all $r<r_0$ and check this for $r=r_0$.

\begin{lem}\label{lem1} It holds \begin{align}
(\p_{\ga_r}...\p_{\ga_1}(\pai
X^A+X^B\Ga^A_{iB}))^\ti&=(\p_{\ga_r}...\p_{\ga_2}((-1)^{(|A|+|B|)|X^B|}R^A_{B\ga_1i}X^B))^\ti,\label{L51a}\end{align}
\begin{multline}(\p_{\ga_r}...\p_{\ga_1}(\p_\ga X^A+(-1)^{|X^B|}X^B\Ga^A_{\ga B}))^\ti=\\ \left\{\begin{array}{ll} 0,&\text{ if }
\ga<\ga_1;\\
(-1)^{s-1}\frac{1}{2}(\p_{\ga_r}...\p_{\ga_{s+1}}\p_{\ga_{s-1}}...\p_{\ga_1}((-1)^{(|A|+|B|)|X^B|}R^A_{B\ga_s\ga_s}X^B))^\ti,&
\text{ if } \ga=\ga_s \text{ for some }s, 1\leq s\leq r;\\
(\p_{\ga_r}...\p_{\ga_2}((-1)^{(|A|+|B|)|X^B|}R^A_{B\ga_1\ga}X^B))^\ti,&\text{ if } \ga_s<\ga<\ga_{s+1} \text{ for some }
s, \\&{ }\qquad\qquad\qquad\qquad 1\leq s\leq r.
\end{array}\right.\label{L51b}
\end{multline}
\end{lem}

{\it Proof.} We have \\ $(\p_{\ga_r}...\p_{\ga_1}(\pai X^A+X^B\Ga^A_{iB}))^\ti$
\begin{align*}=&(\p_{\ga_r}...\p_{\ga_2}(\pai \p_{\ga_1} X^A+(\p_{\ga_1}
X^B)\Ga^A_{iB}+(-1)^{|X^B|}X^B\p_{\ga_1}\Ga^A_{iB}))^\ti\\
=&(\p_{\ga_r}...\p_{\ga_2}(\pai(-(-1)^{|X^B|}X^B\Ga^A_{\ga_1B})-(-1)^{|X^C|}X^C\Ga^B_{\ga_1C}\Ga^A_{iB}+(-1)^{|X^B|}X^B\p_{\ga_1}\Ga^A_{iB}))^\ti\\
=&(\p_{\ga_r}...\p_{\ga_2}(-(-1)^{|X^B|}(\pai X^B)\Ga^A_{\ga_1B}-(-1)^{|X^B|}X^B\pai\Ga^A_{\ga_1B}\\ &
-(-1)^{|X^C|}X^C\Ga^B_{\ga_1C}\Ga^A_{iB}+(-1)^{|X^C|}X^C\p_{\ga_1}\Ga^A_{iC}))^\ti\\
=&(\p_{\ga_r}...\p_{\ga_2}((-1)^{|X^B|}X^C\Ga^B_{iC}\Ga^A_{\ga_1B}-(-1)^{|X^C|}X^C\pai\Ga^A_{\ga_1C}\\ &
-(-1)^{|X^C|}X^C\Ga^B_{\ga_1C}\Ga^A_{iB}+(-1)^{|X^C|}X^C\p_{\ga_1}\Ga^A_{iC}))^\ti\\
=&(\p_{\ga_r}...\p_{\ga_2}((-1)^{|X^C|+|B|+|C|}X^C\Ga^B_{iC}\Ga^A_{\ga_1B}-(-1)^{|X^C|}X^C\pai\Ga^A_{\ga_1C}\\ &
-(-1)^{|X^C|}X^C\Ga^B_{\ga_1C}\Ga^A_{iB}+(-1)^{|X^C|}X^C\p_{\ga_1}\Ga^A_{iC}))^\ti\\
=&(\p_{\ga_r}...\p_{\ga_2}((-1)^{|X^C|}X^CR^A_{C\ga_1i}))^\ti=(\p_{\ga_r}...\p_{\ga_2}((-1)^{|X^C|(|A|+|C|)}R^A_{C\ga_1i}X^C))^\ti\\
=&(\p_{\ga_r}...\p_{\ga_2}((-1)^{(|A|+|B|)|X^B|}R^A_{B\ga_1i}X^B))^\ti.
\end{align*}
Here we used  the induction hypotheses, \eqref{Rcoord} and the fact that  the induction hypotheses imply
$$(\p_{\ga_r}...\p_{\ga_2}(-(-1)^{|X^B|}\pai
X^B))^\ti=(\p_{\ga_r}...\p_{\ga_2}((-1)^{|X^B|}X^C\Ga^B_{iC}))^\ti=(\p_{\ga_r}...\p_{\ga_2}((-1)^{|X^C|+|B|+|C|}X^C\Ga^B_{iC}))^\ti.$$
This proves \eqref{L51a}.

Let us prove \eqref{L51b}. First, if $\ga<\ga_1$, then $(\p_{\ga_r}...\p_{\ga_1}(\p_\ga X^A+(-1)^{|X^B|}X^B\Ga^A_{\ga
B}))^\ti=0$ by the definition of $X^A$. As above, using the induction hypotheses, for any  $t$, $1\leq t\leq r$ we get
\begin{multline}\label{L51pr} (\p_{\ga_r}...\p_{\ga_1}((-1)^{|X^B|}X^B\Ga^A_{\ga B}))^\ti=
(-1)^{t-1}(\p_{\ga_r}...\p_{\ga_{t+1}}\p_{\ga_{t-1}}...\p_{\ga_{1}}\p_{\ga_{t}}((-1)^{|X^B|}X^B\Ga^A_{\ga B}))^\ti\\=
(-1)^{t-1}(\p_{\ga_r}...\p_{\ga_{t+1}}\p_{\ga_{t-1}}...\p_{\ga_{1}}((-1)^{(|A|+|B|)|X^B|}(-(-1)^{|B|+|C|}
\Ga^C_{\ga_tB}\Ga^A_{\ga C}+\p_{\ga_t}\Ga^A_{\ga B})X^B))^\ti.\end{multline} Further, if $\ga=\ga_s$ for some $s$, $1\leq
s\leq r$, then $\p_{\ga_r}...\p_{\ga_1}\p_\ga X^A=0$ and \begin{multline*}(\p_{\ga_r}...\p_{\ga_1}(\p_\ga
X^A+(-1)^{|X^B|}X^B\Ga^A_{\ga B}))^\ti
\\=(-1)^{s-1}(\p_{\ga_r}...\p_{\ga_{s+1}}\p_{\ga_{s-1}}...\p_{\ga_{1}}((-1)^{(|A|+|B|)|X^B|}(-(-1)^{|B|+|C|}
\Ga^C_{\ga_sB}\Ga^A_{\ga_s C}+\p_{\ga_s}\Ga^A_{\ga_s B})X^B))^\ti\\
=(-1)^{s-1}\frac{1}{2}(\p_{\ga_r}...\p_{\ga_{s+1}}\p_{\ga_{s-1}}...\p_{\ga_1}((-1)^{(|A|+|B|)|X^B|}R^A_{B\ga_s\ga_s}X^B))^\ti,\end{multline*}
where we used  \eqref{L51pr} for $t=s$ and \eqref{Rcoord}.

 In the remaining case $\ga_s<\ga<\ga_{s+1}$  for some  $s$,
$1\leq s\leq r$. Then,
\begin{multline*}(\p_{\ga_r}...\p_{\ga_1}\p_\ga X^A)^\ti=(-1)^s(\p_{\ga_r}...\p_{\ga_{s+1}}\p_\ga\p_{\ga_{s}}...\p_{\ga_1}
X^A)^\ti\\=(-1)^s(\p_{\ga_r}...\p_{\ga_{s+1}}\p_\ga\p_{\ga_{s}}...\p_{\ga_2}(-(-1)^{|X^B|}X^B\Ga^A_{\ga_1 B}
))^\ti=(\p_{\ga_r}...\p_{\ga_2}\p_\ga((-1)^{|X^B|}X^B\Ga^A_{\ga_1 B} ))^\ti\\ =(\p_{\ga_r}...\p_{\ga_2}((-1)^{|X^B|}\p_\ga
X^B\Ga^A_{\ga_1 B}+X^B\p_\ga\Ga^A_{\ga_1 B} ))^\ti\\= (\p_{\ga_r}...\p_{\ga_2}(-(-1)^{|X^B|}(-1)^{|X^C|}X^C\Ga^B_{\ga C}
\Ga^A_{\ga_1 B}+X^B\p_\ga\Ga^A_{\ga_1 B} ))^\ti\\=(\p_{\ga_r}...\p_{\ga_2}(-(-1)^{|C|+|B|+(|A|+|B|)|X^B|}\Ga^C_{\ga B}
\Ga^A_{\ga_1 C}X^B+(-1)^{(|A|+|B|)|X^B|}\p_\ga\Ga^A_{\ga_1 B}X^B ))^\ti,\end{multline*} where we used the definition of
$X^A$ and the induction hypotheses. Combining this with \eqref{L51pr} for $t=1$ and \eqref{Rcoord}, we get
$$(\p_{\ga_r}...\p_{\ga_1}(\p_\ga X^A+(-1)^{|X^B|}X^B\Ga^A_{\ga
B}))^\ti=(\p_{\ga_r}...\p_{\ga_2}((-1)^{(|A|+|B|)|X^B|}R^A_{B\ga_1\ga}X^B))^\ti.$$

The lemma is proved. $\Box$

Since $\t\na X_0=0$, for any curve $\gamma$ in $U$ beginning at the point $x\in M$ and ending at a point $y\in U$, we have
$X_y=\tau_\ga X_x$. From this and (ii) it follows that \begin{equation}\label{eq2} (\bar\na^s_{Y_s,...,Y_1}R( Y, Z))^\ti
\t X=0\end{equation} for all $ s\geq 0$ and   $Y,Z,Y_1,...,Y_r\in \T_\M(U)$. Consequently,
\begin{equation}\label{eq3} (\bar\na^{t-1}_{a_t,...,a_2}R^A_{Ba_1a}X^B)^\ti=0\quad\text{for all }  t\geq 1. \end{equation}

Lemma \ref{lem1} and \eqref{eq3} prove \eqref{A1}, \eqref{B1}   and \eqref{st2} for $r=1$. Suppose now that \eqref{st2}
holds for all $r<r_0$ and check it together with \eqref{A1} and \eqref{B1} for $r=r_0$.

\begin{lem}\label{lem2} It holds
$$(\p_{a_r}...\p_{a_{s+1}}((-1)^{(|A|+|B|)|X^B|}\bar\na^{s-1}_{a_s,...,a_2}
R^A_{Ba_1a}X^B))^\ti=(\p_{a_r}...\p_{a_{s+2}}((-1)^{(|A|+|B|)|X^B|}\bar\na^s_{a_{s+1},...,a_2} R^A_{Ba_1a}X^B))^\ti$$ for
all $s$,  $1\leq s\leq r$.
\end{lem}

{\it Proof.} We have \\ $(\p_{a_r}...\p_{a_{s+1}}((-1)^{(|A|+|B|)|X^B|}\bar\na^{s-1}_{a_s,...,a_2} R^A_{Ba_1a}X^B))^\ti$
\begin{align*}
&=(\p_{a_r}...\p_{a_{s+2}}((-1)^{(|A|+|B|)|X^B|}(\p_{a_{s+1}}(\bar\na^{s-1}_{a_s,...,a_2}
R^A_{Ba_1a})X^B\\&\,\,\,\,+(-1)^{|a_{s+1}|(|a_s|+\cdots+|a_1|+|a|+|A|+|B|)}\bar\na^{s-1}_{a_s,...,a_2}
R^A_{Ba_1a}\p_{a_{s+1}}X^B)))^\ti.
\end{align*}
Using \eqref{naRcoord} and the induction hypotheses, we get\\
$(\p_{a_r}...\p_{a_{s+2}}((-1)^{(|A|+|B|)|X^B|}\p_{a_{s+1}}(\bar\na^{s-1}_{a_s,...,a_2} R^A_{Ba_1a})X^B))^\ti$
\begin{align*}
 &=(\p_{a_r}...\p_{a_{s+2}}((-1)^{(|A|+|B|)|X^B|}\bar\na^s_{a_{s+1},...,a_2} R^A_{Ba_1a}X^B\\
 &\,\,\,\,+(-1)^{(|A|+|B|)|X^B|}(-1)^{(|C|+|B|)(|a_{r}|+\cdots|a_1|+|a|)}\Ga^C_{a_{s+1}B}\bar\na^{s-1}_{a_{s},...,a_2}R^A_{Ca_1a}
X^B))^\ti.
\end{align*}
Furthermore,\\
$(\p_{a_r}...\p_{a_{s+2}}((-1)^{(|A|+|B|)|X^B|}(-1)^{|a_{s+1}|(|a_s|+\cdots+|a_1|+|a|+|A|+|B|)}\bar\na^{s-1}_{a_s,...,a_2}
R^A_{Ba_1a}\p_{a_{s+1}}X^B))^\ti$
\begin{align*}
 &=-(\p_{a_r}...\p_{a_{s+2}}((-1)^{(|A|+|B|)|X^B|}(-1)^{|a_{s+1}|(|a_s|+\cdots+|a_1|+|a|+|A|+|B|)}\bar\na^{s-1}_{a_s,...,a_2}
R^A_{Ba_1a}(-1)^{|a_{s+1}||X^C|}X^C\Ga^B_{a_{s+1}C})^\ti\\ &=
-(\p_{a_r}...\p_{a_{s+2}}((-1)^{(|A|+|B|)(|X^C|+|B|+|C|)+|a_{s+1}|(|a_s|+\cdots+|a_1|+|a|+|A|+|B|)+|a_{s+1}||X^C|}\bar\na^{s-1}_{a_s,...,a_2}
R^A_{Ba_1a}X^C\Ga^B_{a_{s+1}C}))^\ti\\ &=
-(\p_{a_r}...\p_{a_{s+2}}((-1)^{(|A|+|C|)(|X^B|+|C|+|B|)+|a_{s+1}|(|a_s|+\cdots+|a_1|+|a|+|A|+|C|)+|a_{s+1}||X^B|}\bar\na^{s-1}_{a_s,...,a_2}
R^A_{Ca_1a}X^B\Ga^C_{a_{s+1}B}))^\ti\\ &= -(\p_{a_r}...\p_{a_{s+2}}((-1)^\sigma\bar\na^{s-1}_{a_s,...,a_2}
R^A_{Ca_1a}X^B\Ga^C_{a_{s+1}B}))^\ti\\ &=
-(\p_{a_r}...\p_{a_{s+2}}((-1)^\sigma(-1)^{(|C|+|a_{s+1}|+|B|)(|a_s|+\cdots+|a_1|+|a|+|C|+|X^B|)}
\Ga^C_{a_{s+1}B}\bar\na^{s-1}_{a_s,...,a_2}
R^A_{Ca_1a}X^B))^\ti,\\&=-(-1)^{(|A|+|B|)|X^B|}(-1)^{(|C|+|B|)(|a_{r}|+\cdots|a_1|+|a|)}\Ga^C_{a_{s+1}B}\bar\na^{s-1}_{a_{s},...,a_2}R^A_{Ca_1a}
X^B))^\ti.
\end{align*} Here we used the fact that by the induction hypotheses $$(\p_{a_r}...\p_{a_{s+2}}(\p_{a_{s+1}}X^B))^\ti=
-(\p_{a_r}...\p_{a_{s+2}}((-1)^{|a_{s+1}||X^C|}X^C\Ga^B_{a_{s+1}C}))^\ti$$ and that
$(-1)^{|a_{s+1}||X^C|}X^C\Ga^B_{a_{s+1}C}=(-1)^{|a_{s+1}|(|X^B|+|B|+|C|)}X^C\Ga^B_{a_{s+1}C}$.

Thus, $$(\p_{a_r}...\p_{a_{s+1}}((-1)^{(|A|+|B|)|X^B|}\bar\na^{s-1}_{a_s,...,a_2}
R^A_{Ba_1a}X^B))^\ti=(\p_{a_r}...\p_{a_{s+2}}((-1)^{(|A|+|B|)|X^B|}\bar\na^s_{a_{s+1},...,a_2} R^A_{Ba_1a}X^B))^\ti.$$ The
lemma is proved. $\Box$

Using Lemma \ref{lem2}, \eqref{eq3} and the fact that if for a function $f\in\O_\M(U)$ holds $|f|=1$, then $\t f=0$, we
get
\begin{multline*}(\p_{a_r}...\p_{a_{s+1}}((-1)^{(|A|+|B|)|X^B|}\bar\na^{s-1}_{a_s,...,a_2}
R^A_{Ba_1a}X^B))^\ti\\=((-1)^{(|A|+|B|)|X^B|}\bar\na^{r-1}_{a_r,...,a_2} R^A_{Ba_1a}X^B)^\ti=(\bar\na^{r-1}_{a_r,...,a_2}
R^A_{Ba_1a}X^B)^\ti=0.\end{multline*} Further,
\begin{multline*}(\p_{\ga_r}...\p_{\ga_1}(\pai
X^A+X^B\Ga^A_{iB}))^\ti=(\p_{\ga_r}...\p_{\ga_2}((-1)^{(|A|+|B|)|X^B|}R^A_{B\ga_1i}X^B))^\ti\\=
((-1)^{(|A|+|B|)|X^B|}\bar\na^{r-1}_{\ga_r,...,\ga_2}R^A_{B\ga_1i}X^B)^\ti=(\bar\na^{r-1}_{\ga_r,...,\ga_2}R^A_{B\ga_1i}X^B)^\ti=0.\end{multline*}
Similarly, $(\p_{\ga_r}...\p_{\ga_1}(\p_\ga X^A+(-1)^{|X^B|}X^B\Ga^A_{\ga B}))^\ti=0$.

 Thus the functions $X^A$
satisfy \eqref{A1} and \eqref{B1}. Consequently, $\na X=0$. From Proposition \ref{prop1} it follows that $X$ does not
depend on the choice of coordinates over $U$. Hence for each coordinate neighbourhood $U\subset M$ we have a unique
parallel section $X\in\E(U)$. Thus we get a parallel section $X\in\E(M)$. The theorem is proved. $\Box$

Recall that the connection $\na$ is called {\it flat} if $\E$ admits local bases of parallel sections.

\begin{corol}
Let $\M=(M,\O_\M)$ be a supermanifold, $\E$ a locally free sheaf of $\O_\M$-supermodules on $\M$ and $\na$ a connection on
$\E$.  Then the following conditions are equivalent: (i) $\na$ is flat; (ii) $R=0$; (iii) $\hol(\na)_x=0$.
\end{corol}

\section{Parallel subsheaves}\label{sps} Let $\M=(M,\O_\M)$ be a supermanifold, $\E$ a locally free sheaf of rank $p|q$ of
$\O_\M$-supermodules on $\M$. For fixed integers $0\leq p_1\leq p$ and $0\leq q_1\leq q$ we assume that
$$\bA,\bB,\bC=1,...,p_1,p+1,...,p+q_1\quad\text{and}\quad \hA,\hB,\hC=p_1+1,...,p,p+q_1+1,...,q.$$ Recall that a subsheaf
$\F\subset\E$ of $\O_\M$-supermodules (of rank $p_1|q_1$)  is called {\it a locally direct subsheaf} if locally there
exists a basis $e_A$ of $\E(U)$ such that $e_\bA$ is a basis of $\F(U)$. For example, a locally direct subsheaf of $\T_\M$
is called {\it a distribution} on $\M$.

Let $\na$ be a connection on $\E$. A locally direct subsheaf $\F\subset\E$ is called {\it parallel} if for any open subset
$U\subset M$ and any $Y\in\T_\M(U)$ and $X\in\F(U)$ it holds $\na_YX\in\F(U)$.

The following theorem is a generalization of Theorem \ref{paralsubbund}.

\begin{theorem}\label{paralsubsh}
Let $\M=(M,\O_\M)$ be a supermanifold, $\E$ a locally free sheaf of $\O_\M$-supermodules on $\M$ and $\na$ a connection on
$\E$.   Then the following conditions are equivalent:
\begin{description} \item[(i)] there exists a  parallel  locally direct subsheaf $\F\subset\E$ of rank $p_1|q_1$;

\item[(ii)] for any point $x\in M$ there exists a vector supersubspace $F_x\subset \E_x$ of dimension $p_1|q_1$
preserved by  $\hol(\na)_x$ and  $\Hol(\t\na)_x$.
\end{description}
\end{theorem}

{\bf Proof of Theorem \ref{paralsubsh}.} Let $\F\subset\E$ be a parallel locally direct subsheaf of rank $p_1|q_1$. Fix a
point $x\in M$. Consider the vector subbundle $F=\cup_{y\in M} \F_y\subset E$. Since $\F\subset\E$ is parallel, the
subbundle $F\subset E$ is parallel. In particular, $F$ is invariant under the parallel displacements in $M$ and for any
curve $\ga$ in $M$ beginning at the point $x\in M$ and ending at a point $y\in M$, we have $F_y=\tau_\ga F_x$. Let
$U\subset M$ be an open subset. Since $F\subset E$ is parallel, for any $X\in\F(U)$ and any $Y,Z,Y_1,...,Y_r\in\T_\M(U)$
we have $$\bar\na^r_{{Y_r},...,{Y_1}}R(Y,Z)X\in\F(U).$$ Consequently, $(\bar\na^r_{{Y_r},...,{Y_1}}R(Y,Z))_y X\in\F_y$ for
all $X\in\F_y$ and $Y,Z,Y_1,...,Y_r\in T_y\M$. This and the fact that $F$ is invariant under the parallel displacements
 prove the implication (i)$\Rightarrow$(ii).

{\it Let us prove the implication} (ii)$\Rightarrow$(i). Suppose that for a point $x\in M$ there exists a vector
supersubspace $F_x\subset \E_x$  preserved by $\hol(\na)_x$ and $\Hol(\t\na)_x$. Then $\Hol(\t\na)_x$ preserves the vector
subspaces $(F_x)_{\bar 0},(F_x)_{\bar 1}\subset E_x$, where $(F_x)_{\bar 0}$ and $(F_x)_{\bar 1}$ are the even and odd
parts of $F_x$, respectively. Consequently, we get parallel vector subbundles $F_{\bar 0},F_{\bar 1}\subset E$ on $M$.
Recall that $E_{\bar 0},E_{\bar 1}\subset E$ are also parallel.  Let $U,x^a$ be a system of local coordinates on $\M$ and
let $e_A$ be a basis of $\Ga(U,E)$ such that $ e_1,..., e_{p_1}$; $ e_{p+1},..., e_{p+q_1}$; $ e_1,..., e_{p}$ and $
e_{p+1},..., e_{p+q}$ are bases of $\Ga(U,F_{\bar 0})$, $\Ga(U,F_{\bar 1})$, $\Ga(U,E_{\bar 0})$ and $\Ga(U,E_{\bar 1})$,
respectively. In particular, $e_A$ is a basis  of $\E(U)$. Since the vector subbundle $F\subset E$ is parallel, we get
\begin{equation}\label{usl0} \t\Ga^\hB_{i\bA}=0.\end{equation} Let $f_\bA^\hB\in\O_\M(U)$ be functions. Consider the
sections $$f_\bA=e_\bA+f_\bA^\hB e_\hB\in\E(U).$$ We will prove  the existence and uniqueness of functions $f_\bA^\hB$
under the conditions $\t f_\bA^\hB=0$, $|f^\hB_\bA|=|\hB|+|\bA|$ and the condition that there exist functions
$X_{a\bA}^\bB\in\O_\M(U)$ such that
\begin{equation}\label{usl}\na_\pa f_\bA=X_{a\bA}^{\bB} f_\bB.\end{equation} Then the supersubmodule
$$\F(U)=\O_\M(U)\otimes\spa_\Real\{f_\bA\}\subset\E(U)$$ will be parallel, i.e. for all $Y\in\T_\M(U)$ and $X\in\F(U)$ it
holds $\na_YX\in\F(U)$. Equation \eqref{usl} is equivalent to the following two equations
\begin{align}\Ga^\bB_{a\bA}+(-1)^{|a|(|\bA|+|\hC|)}f^\hC_\bA\Ga^\bB_{a\hC}&=X^\bB_{a\bA},\label{usl1}\\ \Ga^\hB_{a\bA}+\pa
f^\hB_\bA+(-1)^{|a|(|\bA|+|\hC|)}f^\hC_\bA\Ga^\hB_{a\hC}&=X^\bD_{a\bA}f^\hB_\bD.\label{usl2}\end{align} Combining
\eqref{usl1} and \eqref{usl2}, we get
\begin{equation}\label{AABB}
\Ga^\hB_{a\bA}+\pa f^\hB_\bA+(-1)^{|a|(|\bA|+|\hC|)}f^\hC_\bA\Ga^\hB_{a\hC}-(\Ga^\bD_{a\bA}+
(-1)^{|a|(|\bA|+|\hC|)}f^\hC_\bA\Ga^\bD_{a\hC})f^\hB_\bD=0.\end{equation} We will show that Equation \eqref{AABB} has a
unique solution $f_\bA^\hB$ satisfying the conditions $\t f_\bA^\hB=0$ and $|f^\hB_\bA|=|\hB|+|\bA|$, then \eqref{usl}
will hold for the functions $X^\bB_{a\bA}$ given by \eqref{usl1}. Equation \eqref{AABB} can be written as
\begin{align}(\p_{\ga_r}...\p_{\ga_1}(\Ga^\hB_{i\bA}+\pai f^\hB_\bA+f^\hC_\bA\Ga^\hB_{i\hC}-(\Ga^\bD_{i\bA}+
f^\hC_\bA\Ga^\bD_{i\hC})f^\hB_\bD))^\ti&=0,\label{AA}\\ (\p_{\ga_r}...\p_{\ga_1}(\Ga^\hB_{\ga\bA}+\p_\ga
f^\hB_\bA+(-1)^{|\bA|+|\hC|}f^\hC_\bA\Ga^\hB_{\ga\hC}-(\Ga^\bD_{\ga\bA}+
(-1)^{|\bA|+|\hC|}f^\hC_\bA\Ga^\bD_{\ga\hC})f^\hB_\bD))^\ti&=0,\label{BB}\end{align} where $0\leq r\leq m$. The further
proof is similar to the proof of Theorem \ref{paralsec}. As in Section \ref{parallel}, we can use \eqref{BB} to define
functions  $f_{\bA\ga\ga_1...\ga_r}^\hB$ ($\ga<\ga_1<\cdots<\ga_r$, $0\leq r\leq m-1$) and $f_\bA^\hB$ using the condition
$\t f_\bA^\hB=0$. Then these functions satisfy $|f^\hB_\bA|=|\hB|+|\bA|$. We must prove that $f_\bA^\hB$ satisfy
\eqref{AA} and \eqref{BB} for $\ga_1<\cdots<\ga_r$, $0\leq r\leq m$ and all $\ga$. We will prove this by the induction
over $r$. Parallel to this we will prove that
\begin{multline}\label{stst2} (\p_{a_r}...\p_{a_{s+1}}(\bar\na^{s-1}_{a_s,...,a_2}
R^\hB_{\bA a_1a}+(-1)^{|\bA|+|\hC|}f^\hC_\bA \bar\na^{s-1}_{a_s,...,a_2}R^\hB_{\hC
a_1a}\\-\bar\na^{s-1}_{a_s,...,a_2}R^\bD_{\bA a_1a}f^\hB_\bD-(-1)^{|\bA|+|\hC|}f^\hC_\bA
\bar\na^{s-1}_{a_s,...,a_2}R^\hD_{\hC a_1a}f^\hB_\bD))^\ti=0
\end{multline}for all  $r\geq 1$  and   $1\leq s\leq r$.
For $r=0$  Equation \eqref{AA} follows from \eqref{usl0} and the condition $\t f_\bA^\hB=0$; Equation \eqref{BB} follows
from the definition of the functions $f_{\bA\ga}^\hB$.
\begin{lem}\label{lemAA1} For $r\geq 1$ it holds \begin{multline*}(\p_{\ga_r}...\p_{\ga_1}
(\Ga^\hB_{i\bA}+\pai f^\hB_\bA+f^\hC_\bA\Ga^\hB_{i\hC}-(\Ga^\bD_{i\bA}+
f^\hC_\bA\Ga^\bD_{i\hC})f^\hB_\bD))^\ti\\=(\p_{\ga_r}...\p_{\ga_2}( R^\hB_{\bA \ga_1i}+(-1)^{|\bA|+|\hC|}f^\hC_\bA
R^\hB_{\hC \ga_1i}-R^\bD_{\bA \ga_1i}f^\hB_\bD-(-1)^{|\bA|+|\hC|}f^\hC_\bA R^\hD_{\hC
\ga_1i}f^\hB_\bD))^\ti.\end{multline*} A similar holds for $(\p_{\ga_r}...\p_{\ga_1}(\Ga^\hB_{\ga\bA}+\p_\ga
f^\hB_\bA+(-1)^{|\bA|+|\hC|}f^\hC_\bA\Ga^\hB_{\ga\hC}-(\Ga^\bD_{\ga\bA}+
(-1)^{|\bA|+|\hC|}f^\hC_\bA\Ga^\bD_{\ga\hC})f^\hB_\bD))^\ti$ (see Lemma \ref{lem1}). \end{lem}

{\it The proof} is similar to the proof of Lemma \ref{lem1}. $\Box$

\begin{lem}\label{lemAA2} For $r\geq 1$ it holds
$$(\bar\na^{r-1}_{a_r,...,a_2} R^\hB_{\bA a_1a})^\ti=0.$$\end{lem}

{\it The proof} follows from the facts that the holonomy algebra $\hol(\na)_x$ preserves the vector subspace $F_x\subset
\E_x$ and that the distribution $F\subset E$ is parallel along all curves in $M$. $\Box$

Lemma \ref{lemAA2} proves \eqref{AA}, \eqref{BB} and \eqref{stst2} for $r=1$. Fix $r_0\geq 2$. Suppose that \eqref{AA},
\eqref{BB} and \eqref{stst2} hold for all $r<r_0$ and check this for $r=r_0$.

\begin{lem}\label{lemAA3}It holds \begin{align*}
&(\p_{a_r}...\p_{a_{s+1}}(\bar\na^{s-1}_{a_s,...,a_2} R^\hB_{\bA a_1a}+(-1)^{|\bA|+|\hC|}f^\hC_\bA
\bar\na^{s-1}_{a_s,...,a_2}R^\hB_{\hC a_1a}\\&\qquad-\bar\na^{s-1}_{a_s,...,a_2}R^\bD_{\bA
a_1a}f^\hB_\bD-(-1)^{|\bA|+|\hC|}f^\hC_\bA \bar\na^{s-1}_{a_s,...,a_2}R^\hD_{\hC a_1a}f^\hB_\bD))^\ti\\
&=(\p_{a_r}...\p_{a_{s+2}}(\bar\na^{s}_{a_{s+1},...,a_2} R^\hB_{\bA a_1a}+(-1)^{|\bA|+|\hC|}f^\hC_\bA
\bar\na^{s}_{a_{s+1},...,a_2}R^\hB_{\hC a_1a}\\&\qquad-\bar\na^{s}_{a_{s+1},...,a_2}R^\bD_{\bA
a_1a}f^\hB_\bD-(-1)^{|\bA|+|\hC|}f^\hC_\bA \bar\na^{s}_{a_{s+1},...,a_2}R^\hD_{\hC a_1a}f^\hB_\bD))^\ti\end{align*} for
all $s$, $1\leq s\leq r$.
\end{lem}

{\it The proof} is similar to the proof of Lemma \ref{lem2}. $\Box$

Now Equations \eqref{AA}, \eqref{BB} and \eqref{stst2} follow from the above lemmas and the induction hypotheses.

We have proved that the supersubmodule $\F(U)=\O_\M(U)\otimes\spa_\Real\{f_\bA\}\subset\E(U)$ is parallel. We claim that
$\F(U)$ does not depend on the choice of the basis $e_A$. Suppose that we have another basis $e'_A$ of $\Ga(U,E)$ with the
same property as above. Then there exist functions $H^\bB_\bA,H^B_\hA\in\O_M(U)$ such that $e'_\bA=H^\bB_\bA e_\bB$ and
$e'_\hA=H^B_\hA e_B$. Furthermore, $$f'_\bA=e'_\bA+{f'}_\bA^\hB e'_\hB=H^\bB_\bA e_\bB+{f'}_\bA^\hB H^B_\hB e_B=
(H^\bB_\bA+{f'}_\bA^\hC H^\bB_\hC) e_\bB+{f'}_\bA^\hC H^\hB_\hC e_\hB.$$ Since $(H^\bB_\bA+{f'}_\bA^\hC
H^\bB_\hC)^\ti=H^\bB_\bA$ is an invertible matrix, the matrix $H^\bB_\bA+{f'}_\bA^\hC H^\bB_\hC$ is also invertible.
Denote by $Y^\bB_\bA$ its inverse matrix. Then, $Y^\bA_\bC f'_\bA=e_\bC+Y^\bA_\bC {f'}_\bA^\hC H^\hB_\hC e_\hB.$ Moreover,
$$\na_\pa (Y^\bA_\bC f'_\bA)=(\pa Y^\bB_\bC+(-1)^{|a||Y^\bA_\bC|}Y^\bA_\bC {X'}^\bB_{a\bA})f'_\bB=(\pa
Y^\bB_\bC+(-1)^{|a||Y^\bA_\bC|}Y^\bA_\bC {X'}^\bB_{a\bA})(H^\bD_\bB+{f'}_\bB^\hC H^\bD_\hC)(Y^\bA_\bD f'_\bA).$$ From the
uniqueness of the functions $f^\hB_\bC$ it follows that $f^\hB_\bC=Y^\bA_\bC {f'}_\bA^\hC H^\hB_\hC$. Consequently,
$f_\bC=Y^\bA_\bC f'_\bA.$ Similarly, there exist functions ${Y'}^\bA_\bC$ such that $f'_\bC={Y'}^\bA_\bC f_\bA.$ This
proves that $\F(U)$ does not depend on the choice of the basis $e_A$. Note that the sections $f_\bA,e_\hA$ form a basis of
$\E(U)$. Thus we get a parallel locally direct subsheaf $\F\subset\E$ of rank $p_1|q_1$. The theorem is proved. $\Box$

\section{Holonomy of linear connections over supermanifolds}\label{LinConnect}

Let $\M=(M,\O_\M)$ be a supermanifold of dimension $n|m$. In this section we consider a connection $\na$ on the tangent
sheaf $\T_\M$ of $\M$. Then in Definition \ref{defhol} of the holonomy algebra $\hol(\na)_x$ we may choose $\bar\na=\na$.
If we put $\E=\T_\M$, then in the above notation, $E=\cup_{y\in M}T_y\M=T\M$. In particular, $E_{\bar 0}=TM$ is the
tangent bundle over $M$. We get the connections $\tilde\na$ and $\tilde\na|_{TM}$ on the vector bundles $T\M$ and $TM$,
respectively. We identify the holonomy algebra $\hol(\na)_x$ and the group $\Hol(\tilde\na)_x$ with a supersubalgebra
$\hol(\na)\subset\gl(n|m,\Real)$ and a Lie subgroup  $\Hol(\tilde\na)\subset \GL(n,\Real)\times \GL(m,\Real)$,
respectively.

On the cotangent sheaf $\T^*_\M=\Hom_{\O_\M}(\T_\M,\O_\M)$ of $\M$ we get the connection $\na^*$ defined as follows
$$(\na^*_X\varphi)Y=X\varphi(Y)-(-1)^{|X||\varphi|}\varphi(\na_X Y),$$ where $U\subset M$ is an open subset, $X,Y\in
\T_\M(U)$ and $\varphi\in\T^*_\M(U)$ are homogeneous.  The curvature tensor $R^*$ of the connection $\na^*$ is given by
$$(R^*(Y,Z)\varphi)X=-(-1)^{|\varphi|(|Y|+|Z|)}\varphi(R(Y,Z)X),$$ where $U\subset M$ is an open subset, $X,Y\in \T_\M(U)$
and $\varphi\in\T^*_\M(U)$ are homogeneous. Let $\ga:[a,b]\subset\Real\to M$ be a curve, $\tau_\ga:T_{\ga(a)}\M\to
T_{\ga(b)}\M$   and $\tau^*_\ga:T^*_{\ga(a)}\M\to T^*_{\ga(b)}\M$ be the parallel displacements of the connections $\na$
and $\na^*$, respectively. Then, $$(\tau^*_\ga\varphi)X=\varphi(\tau^{-1}_\ga X),$$ where $\varphi\in\T^*_{\ga(a)}\M$ and
$X\in \T_{\ga(b)}\M$.

Consider the sheaf of tensor fields of type $(r,s)$ over $\M$,
$$\T^{r,s}_\M=\otimes^r_{\O_\M}\T_{\M}\bigotimes_{\O_\M}\otimes^s_{\O_\M}\T^*_{\M}.$$ Define the connection $\na^{r,s}$ on
this sheaf by
\begin{multline*}\na^{r,s}_X(X_1\otimes\cdots\otimes X_r\otimes \varphi_1\otimes\cdots\otimes
\varphi_s)\\=\sum_{i=1}^{r}(-1)^{|X|(|X_1|+\cdots+|X_{i-1}|)}X_1\otimes\cdots\otimes X_{i-1}\otimes\na_XX_i \otimes
X_{i+1} \otimes\cdots\otimes X_r\otimes \varphi_1\otimes\cdots\otimes \varphi_s\\+
\sum_{j=1}^{s}(-1)^{|X|(|X_1|+\cdots+|X_{r}|+|\varphi_1|+\cdots+|\varphi_{j-1}|)}X_1\otimes\cdots\otimes
X_r\otimes\varphi_1\otimes\cdots\otimes \varphi_{j-1}\otimes\na^*_X\varphi_j \otimes\varphi_{j+1}
\otimes\cdots\otimes\varphi_s,\end{multline*} where $U\subset M$ is an open subset , $X,X_1,...,X_r\in \T_\M(U)$ and
$\varphi_1,...,\varphi_r\in\T^*_\M(U)$ are homogeneous. For the curvature tensor $R^{r,s}$ of this connection we get
\begin{multline*}R^{r,s}(Y,Z)(X_1\otimes\cdots\otimes X_r\otimes \varphi_1\otimes\cdots\otimes
\varphi_s)\\=\sum_{i=1}^{r}(-1)^{(|Y|+|Z|)(|X_1|+\cdots+|X_{i-1}|)}X_1\otimes\cdots\otimes X_{i-1}\otimes R(Y,Z)X_i
\otimes X_{i+1} \otimes\cdots\otimes X_r\otimes \varphi_1\otimes\cdots\otimes
\varphi_s\\+\sum_{j=1}^{s}(-1)^{(|Y|+|Z|)(|X_1|+\cdots+|X_{r}|+|\varphi_1|+\cdots+|\varphi_{j-1}|)}X_1\otimes\cdots\otimes
X_r\otimes\varphi_1\otimes\cdots\otimes \varphi_{j-1}\otimes R^*(Y,Z)\varphi_j \otimes\varphi_{j+1}
\otimes\cdots\otimes\varphi_s,\end{multline*} where $U\subset M$ is an open subset, $Y,Z,X_1,...,X_r\in \T_\M(U)$ and
$\varphi_1,...,\varphi_s\in\T^*_\M(U)$ are homogeneous. For the parallel displacement $\tau^{r,s}_\ga$ of the connections
$\na^{r,s}$ along a curve $\ga:[a,b]\subset\Real\to M$ it holds $\tau^{r,s}_\ga(X_1\otimes\cdots\otimes X_r\otimes
\varphi_1\otimes\cdots\otimes \varphi_s)=(\tau_\ga X_1\otimes\cdots\otimes\tau_\ga X_r\otimes\tau^*_\ga
\varphi_1\otimes\cdots\otimes \tau^*_\ga\varphi_s),$ where $X_1,...,X_r\in T_{\ga(a)}\M$ and $\varphi_1,...,\varphi_s\in
T^*_{\ga(a)}\M$.

Thus we see that the holonomy algebra $\hol(\na^{r,s})_x$ of the connection $\na^{r,s}$ at a point $x\in M$ and the group
$\Hol(\tilde\na^{r,s})_x$ coincide with the tensor extension of the representation of the holonomy algebra $\hol(\na)_x$
and with the tensor extension of the representation of the group $\Hol(\tilde\na)_x$, respectively. From Theorem
\ref{paralsec} we immediately read the following.

\begin{theorem}\label{paraltens}
Let $\M=(M,\O_\M)$ be a supermanifold and $\na$ a connection on $\T_\M$.   Then the following conditions are equivalent:
\begin{description} \item[(i)] there exists a non-zero paralel tensor $P\in\T^{r,s}_\M(M)$;

\item[(ii)] for any point $x\in M$ there exists a non-zero tensor $P_x\in T^{r,s}_x\M$ annihilated by the tensor extension of the representation of the holonomy
algebra $\hol(\na)_x$ and preserved by tensor extension of the representation of the group $\Hol(\t\na)_x$.
\end{description}
%If the manifold  $M$ is simply connected, then  condition (ii) can be written as follows
%\begin{description}
%\item[(ii)] for any point $x\in M$ there exists a non-zero vector $X_x\in \E_x$ annihilates by the holonomy algebra $\hol(\na)_x$.
%\end{description}
\end{theorem}

\begin{ex} In Table \ref{tab1} we give equivalent conditions for existance of some parallel tensors on supermanifolds and
inclusions of the holonomy. The supersubalgebras of $\gl(n|m,\Real)$ that appear in the table are exactly the
supersubalgebras annihilating the corresponding tensor at one point.
\end{ex}

\begin{tab}\label{tab1} Examples of parallel structures and the corresponding holonomy
\begin{longtable}{|l|c|c|c|}\hline parallel structure on $\M$& $\hol(\na)$ is&
$\Hol(\tilde\na)$ is&restriction\\&contained in&contained  in&\\\hline

Riemannian supermetric,&$\osp(p_0,q_0|2k)$&$\Ort(p_0,q_0)\times \Sp(2k,\Real)$&$n=p_0+q_0,m=2k$\\ i.e. even
non-degenerate&&&\\ supersymmetric metric&&&\\\hline

 even non-degenerate &$\osp^{\rm sk}(2k|p,q)$&$\Sp(2k,\Real)\times \Ort(p,q)$&$n=2k,m=p+q$\\  super skew-symmetric metric&&&\\\hline

complex structure&$\gl(k|l,\Co)$&$\GL(k,\Co)\times \GL(l,\Co)$&$n=2k,l=2m$\\\hline

odd complex structure,&$\q(n,\Real)$&$\big\{\big(\begin{smallmatrix} A&0\\0&A\end{smallmatrix}\big)\big| A\in
\GL(n,\Real)\big\}$&$m=n$\\ i.e. odd automorphism &&&\\$J$ of $\T_\M$ with $J^2=-\id$ &&& \\\hline

odd non-degenerate &$\pe(n,\Real)$&$\big\{\big(\begin{smallmatrix} A&0\\0&A\end{smallmatrix}\big)\big| A\in
\GL(n,\Real)\big\}$&$m=n$\\supersymmetric metric &&&\\\hline

odd non-degenerate super&$\pe^{sk}(n,\Real)$&$\big\{\big(\begin{smallmatrix} A&0\\0&A\end{smallmatrix}\big)\big| A\in
\GL(n,\Real)\big\}$&$m=n$\\ skew-symmetric metric &&&\\\hline

\end{longtable}
\end{tab}

{\it The torsion} of the connection $\na$ is given by the formula \begin{equation}\label{torsion}
T(X,Y)=\na_XY-(-1)^{|X||Y|}\na_YX-[X,Y],\end{equation} where  $U\subset M$ is open and $X,Y\in \T_\M(U)$ are homogeneous.
If the connection $\na$ is torsion-free, i.e. $T=0$, then the curvature tensor $R$ satisfies the Bianchi identity
\begin{equation}\label{Bianchi} R(X,Y)Z+(-1)^{|X|(|Y|+|Z|)}R(Y,Z)X+(-1)^{|Z|(|X|+|Y|)}R(Z,X)Y=0\end{equation} for  $U\subset M$  open and all homogeneous $X,Y,Z\in \T_\M(U)$.

Recall that {\it a distribution} on $\M$ is a locally direct subsheaf $\F\subset\T_\M$. The distribution $\F$ is called
{\it involutive}, if $[\F(U),\F(U)]\subset\F(U)$ for all open subsets $U\subset M$. Theorem \ref{paralsubsh}
 gives us a one-to-one
correspondence between parallel distributions on $\M$ and vector supersubspaces $F\subset T_x\M$ preserved by
$\hol(\na)_x$ and $\Hol(\tilde\na)_x$. From \eqref{torsion} it follows that if $\na$ is torsion-free, then any parallel
distribution on $\M$ is involutive.

\section{Berger superalgebras}\label{secBerger}

Let $V$ be a real or complex vector superspace and  $\g\subset\gl(V)$ a supersubalgebra. {\it The space of algebraic
curvature tensors of type} $\g$ is the vector superspace $\R(\g)=\R(\g)_{\bar 0}\oplus\R(\g)_{\bar 1},$ where
$$\R(\g)=\left\{R\in V^*\wedge
V^*\otimes\g\left|\begin{matrix}R(X,Y)Z+(-1)^{|X|(|Y|+|Z|)}R(Y,Z)X+(-1)^{|Z|(|X|+|Y|)}R(Z,X)Y=0\\\text{for all homogeneous
} X,Y,Z\in V
\end{matrix}\right\}\right..$$
Obviously, $\R(\g)$ is a $\g$-module with respect to the action \begin{equation}\label{R_A}A\cdot R=R_A,\quad
R_A(X,Y)=[A,R(X,Y)]-(-1)^{|A||R|}R(AX,Y)-(-1)^{|A|(|R|+|X|)}R(X,AY),\end{equation} where $A\in\g$, $R\in\R(\g)$ and
$X,Y\in V$ are homogeneous. If $\M$ is a supermanifold and $\na$ is a linear torsion-free connection on $\T_\M$, then
applying covariant derivatives to  \eqref{Bianchi}, we get  that $(\na^r_{Y_r,...,Y_1}R)_x\in\R(\hol(\na)_x)$ for all
$r\geq 0$ and $Y_1,...,Y_r\in T_xM$. Moreover, $|(\na^r_{Y_r,...,Y_1}R)_x|=|Y_1|+\cdots+|Y_r|$, whenever $Y_1,...,Y_r$ are
homogeneous.

Define the vector supersubspace $$L(\R(\g))=\spa\{R(X,Y)|R\in\R(\g),\,\,X,Y\in V\}\subset \g.$$ From \eqref{R_A} it
follows that $L(\R(\g))$ is an ideal in $\g$. We call a supersubalgebra $\g\subset\gl(V)$ {\it a Berger superalgebra} if
$L(\R(\g))=\g$.

\begin{prop} Let $\M$ be a supermanifold of dimension $n|m$ with a linear torsian-free connection $\na$.
Then its holonomy algebra $\hol(\na)\subset\gl(n|m,\Real)$ is a Berger superalgebra.\end{prop}

{\bf Proof.} The proof follows from Definition \ref{defhol} and \eqref{Bianchi}. $\Box$

Consider the vector superspace $$\R^\na(\g)= \left\{S\in V^*
\otimes\R(g)\left|\begin{matrix}S_X(Y,Z)+(-1)^{|X|(|Y|+|Z|)}S_Y(Z,X)+(-1)^{|Z|(|X|+|Y|)}S_Z(X,Y)=0\\ \text{for all
homogeneous } X,Y,Z\in V
\end{matrix}\right\}\right..$$ If $\M$ is a supermanifold and $\na$ is a linear torsion-free connection on $\T_\M$, then
applying covariant derivatives to  the second Bianchi identity, we get  that
$(\na^r_{Y_r,...,Y_2,\cdot}R)_x\in\R^\na(\hol(\na)_x)$ for all $r\geq 1$ and $Y_2,...,Y_r\in T_xM$. Moreover,
$|(\na^r_{Y_r,...,Y_2,\cdot}R)_x|=|Y_2|+\cdots+|Y_r|$, whenever $Y_2,...,Y_r$ are homogeneous.

\section{Holonomy of locally symmetric superspaces}\label{secsym}

Let $\M$ be a supermanifold of dimension $n|m$ with a linear connection $\na$. The supermanifold $(\M,\na)$ is called {\it
locally symmetric} if $\na$ is torsion-free and $\na R=0$. Note that in this situation the underlying manifold
$(M,\t\na|_{TM})$ is locally symmetric as well. Let $x\in M$. Theorem \ref{paraltens} implies that $\hol(\na)_x$
annihilates the value $R_x\in\R(\hol(\na)_x)$ and $\Hol(\t\na)_x$ preserves $R_x$. We get that
$\hol(\na)_x=\spa\{R_x(X,Y)|X,Y\in T_x\M\}.$

More generally, let $V$ be a vector superspace and suppose that $\g\subset\gl(V)$ is a subalgebra that annihilates an
$R\in\R(\g)$. Consider the Lie superalgebra $\h=\g+V$ with the Lie brackets $$[x,y]=-R(x,y),\quad [A,x]=Ax,\quad
[A,B]=A\circ B-(-1)^{|A||B|}B\circ A,$$ where $x,y\in V$ and $A,B\in\g$. We get that $\h=\g+V$ is {\it a symmetric
decomposition} of the Lie superalgebra $\h$ \cite{Cortes2}. All such decompositions for $\h$ simple are described in
\cite{Serganova}. In particular, this allows to find all possible irreducible holonomy algebras of locally symmetric
supermanifolds.

A Berger superalgebra $\g$ is called {\it symmetric} if $\R^\na(\g)=0$.

\begin{prop} Let $\M$ be a  supermanifold with a torsian free connection $\na$. If $\hol(\na)$ is a symmetric Berger
superalgebra, then $(\M,\na)$ is locally symmetric.\end{prop}

{\bf Proof.} We need to prove that $\na R=0$.  Since $\R^\na(\hol(\na)_y)=0$ for all $y\in M$, we get
$(\na^r_{Y_r,...,Y_1}R)^\ti$ for all $r\geq 0$ and all vector fields $Y_1,...,Y_r$ on $\M$. Using this, \eqref{naRcoord}
and double induction over $r$ and $s$  it is easy to get that
$(\p_{a_r}\cdots\p_{a_{s+1}}\na^{s+1}_{a_s,...,a_1,f}R^d_{cab})^\ti=0$ for all $r$ and $s$ such that $r\geq s\geq 0$. In
particular, $(\p_{\ga_r}\cdots\p_{\ga_{1}}\na_{f}R^d_{cab})^\ti=0$ for all $r\geq 0$. Thus, $\na R=0$. $\Box$

The proof of the following proposition is as in \cite{Sch}.

\begin{prop} Let $\g\subset\gl(V)$ be an irreducible Berger superalgebra. If $\g$ annihilates the module $\R(\g)$,
then $\g$ is a symmetric Berger superalgebra.
\end{prop}

\section{Holonomy of Riemannian supermanifolds} \label{holRiem}

{\it A Riemannian supermanifold} $(\M,g)$ is a supermanifold $\M$ of dimension $n|m$, $m=2k$ endowed with an even
non-degenerate supersymmetric metric $g$ \cite{Cortes1}. In particular, the value $g_x$ of $g$ at a point $x\in M$
satisfies: $g_x|_{(T_x\M)_{\bar 0},(T_x\M)_{\bar 1}}=0$, $g_x|_{(T_x\M)_{\bar 0}\times (T_x\M)_{\bar 0}}$ is
non-degenerate, symmetric and $g_x|_{(T_x\M)_{\bar 1}\times (T_x\M)_{\bar 1}}$ is non-degenerate, skew-symmetric. The
metric $g$ defines a pseudo-Riemannian metric $\tilde g$ on the manifold $M$. Note that $\tilde g$ is not assumed to be
positively defined. The supermanifold $(\M,g)$ has a unique linear connection  $\na$ such that $\na$ is torsion-free and
$\na g=0$. This connection is called {\it the Levi-Civita connection}. We denote the holonomy algebra of the connection
$\na$ by $\hol(\M,g)$. As we have already noted, $\hol(\M,g)\subset\osp(p_0,q_0|2k)$ and $\Hol(\tilde\na)\subset
\Ort(p_0,q_0)\times \Sp(2k,\Real)$, where $(p_0,q_0)$ is the signature of the pseudo-Riemannian metric $\tilde g$.

The  K\"ahlerian, hyper-K\"ahlerian and quaternionic-K\"ahlerian supermanifolds are defined in the natural way, see e.g
\cite{Cortes1}. We define {\it special K\"ahlerian} or {\it Calabi-Yau} supermanifolds by the condition from Table
\ref{tab2}.

\begin{ex} In Table \ref{tab2} we give equivalent conditions for special geometry of $(\M,g)$ and
inclusion of the holonomy.
%The definitions of supersubalgebras of $\osp(n|2k,\Real)$ that appear in the table are clear from
%the context.
\end{ex}

\begin{tab}\label{tab2} Special geometries of Riemannian supermanifolds and the corresponding holonomies
\begin{longtable}{|l|c|c|c|}\hline type of $(\M,g)$& $\hol(\M,g)$ is&
$\Hol(\tilde\na)$ is&restriction\\&contained in&contained  in&\\\hline

K\"ahlerian &$\u(p_0,q_0|p_1,q_1)$& $\Un(p_0,q_0)\times \Un(p_1,q_1)$&$n=2p_0+2q_0,$\\ &&&$m=2p_1+2q_1$\\\hline

special K\"ahlerian &$\su(p_0,q_0|p_1,q_1)$& $\Co^*(\SU(p_0,q_0)\times \SU(p_1,q_1))$&$n=2p_0+2q_0$,\\
&&&$m=2p_1+2q_1$\\\hline

hyper-K\"ahlerian &$\hosp(p_0,q_0|p_1,q_1)$& $\Sp(p_0,q_0)\times \Sp(p_1,q_1)$&$n=4p_0+4q_0$,\\ &&& $m=4p_1+4q_1$\\\hline

quaternionic- &$\sp(1)\oplus\hosp(p_0,q_0|p_1,q_1)$& $\Sp(1)(\Sp(p_0,q_0)\times \Sp(p_1,q_1))$&$n=4p_0+4q_0\geq 8$,\\
K\"ahlerian&&& $m=4p_1+4q_1$\\\hline

\end{longtable}
\end{tab}

 Define the Ricci tensor $\Ric$ of the supermanifold $(\M,g)$ by the formula \begin{equation}\label{defRic}\Ric(Y,Z)=\str(X\mapsto
(-1)^{|X||Z|}R(Y,X)Z),\end{equation} where $U\subset M$ is open and $X,Y,Z\in\T_\M(U)$ are homogeneous.

\begin{prop} Let $(\M,g)$ be a  K\"ahlerian supermanifold, then $\Ric=0$ if and only if $\hol(\M,g)\subset\su(p_0,q_0|p_1,q_1)$.
In particular, if $(\M,g)$ is special  K\"ahlerian, then $\Ric=0$; if $M$ is simply connected, $(\M,g)$ is  K\"ahlerian
and $\Ric=0$, then $(\M,g)$ is special  K\"ahlerian. \end{prop}

{\bf Proof.} If $R$ is an algebraic curvature tensor of type $\u(p_0,q_0|p_1,q_1)$, then $\Ric(R)$ is defined by
\eqref{defRic}. The following formula can be proved as in the usual case up to additional signs
\begin{equation}\label{propRic}
\Ric(R)(Y,Z)=\frac{1}{2}\str(J\circ R(JY,Z)),\end{equation} where $J$ is the complex structure. Recall that
\begin{equation}\label{su}\su(p_0,q_0|p_1,q_1)=\{\xi\in\u(p_0,q_0|p_1,q_1)|\str(J\circ \xi)=0\}.\end{equation}
Combining  \eqref{naR} and \eqref{naRcoord} and using the equality \begin{multline*}
(-1)^{|a_r|(|a_{r-1}|+\cdots+|a_1|+|a|+|b|+|d|+|c|)+|b|(|d|+1)} \na^{r-1}_{a_{r-1},...,a_1}R^c_{dab}\Ga^b_{a_rc}\\-
(-1)^{(|c|+|b|)(|a_{r-1}|+\cdots+|a_1|+|a|)+|b|(|d|+1)}\Ga^c_{a_rb}\bar\na^{r-1}_{a_{r-1},...,a_1}R^b_{dac}=0,\end{multline*}
we get
\begin{align}\Ric(\na^s_{\p_{a_s},...,\p_{a_1}}R)(\p_a,\p_d)=
&\p_{a_s}(\Ric(\na^{s-1}_{\p_{a_{s-1}},...,\p_{a_1}}R)(\p_a,\p_d))-\Ric(\na^{s-1}_{\na_{\p_{a_s}}\p_{a_{s-1}},...,\p_{a_1}}R)(\p_a,\p_d)\nonumber\\
&-(-1)^{|a_s||a_{s-1}|}\Ric(\na^{s-1}_{\p_{a_{s-1}},\na_{\p_{a_s}}\p_{a_{s-2}},...,\p_{a_1}}R)(\p_a,\p_d)\nonumber
\\&-\cdots-(-1)^{|a_s|(|a_{s-1}|+\cdots+|a_2|)}
\Ric(\na^{s-1}_{\p_{a_{s-1}},...,\p_{a_2},\na_{\p_{a_s}}\p_{a_1}}R)(\p_a,\p_d)\label{naRic}\\
&-(-1)^{|a_s|(|a_{s-1}|+\cdots+|a_1|)}\Ric(\na^{s-1}_{\p_{a_{s-1}},...,\p_{a_1}}R)(\na_{\p_{a_s}}\p_a,\p_d)\nonumber\\
&-(-1)^{|a_s|(|a_{s-1}|+\cdots+|a_1|+|a|)}\Ric(\na^{s-1}_{\p_{a_{s-1}},...,\p_{a_1}}R)(\p_a,\na_{\p_{a_s}}\p_d).\nonumber
\end{align}
Suppose that $\Ric=0$. Using \eqref{naRic} and induction it is easy to prove that
$\Ric(\na^s_{\p_{a_s},...,\p_{a_1}}R)(\p_a,\p_d)=0$ for all $s\geq 0$. Consequently,
$(\Ric(\na^s_{\p_{a_s},...,\p_{a_1}}R)(\p_a,\p_d))^\ti=0$ for all $s\geq 0$. This, \eqref{propRic} and \eqref{su} yield
$\hol(\M,g)\subset\su(p_0,q_0|p_1,q_1)$.

Suppose that $\hol(\M,g)\subset\su(p_0,q_0|p_1,q_1)$. Then $(\Ric(\na^s_{\p_{a_s},...,\p_{a_1}}R)(\p_a,\p_d))^\ti=0$ for
all $s\geq 0$.  Using this, \eqref{naRic} and double induction over $r$ and $s$, it is easy to show that
$$(\p_{a_r}\cdots\p_{a_{s+1}}(\Ric(\na^s_{\p_{a_s},...,\p_{a_1}}R)(\p_a,\p_d)))^\ti=0.$$ for all $r$ and $s$  such that
$r\geq s\geq 0$. In particular, $(\p_{\ga_r}\cdots\p_{\ga_{1}}(\Ric(\p_a,\p_d)))^\ti=0$ for all $r\geq 0$, i.e. $\Ric=0$.
The proposition is proved. $\Box$

\begin{corol} Let $(\M,g)$ be a  quaternionic-K\"ahlerian supermanifold, then $\Ric=0$ if and only if
$\hol(\M,g)\subset\hosp(p_0,q_0|p_1,q_1)$. In particular, if $(\M,g)$ is hyper-K\"ahlerian, then $\Ric=0$; if $M$ is
simply connected, $(\M,g)$ is quaternionic-K\"ahlerian and $\Ric=0$, then $(\M,g)$ is hyper-K\"ahlerian. \end{corol}

We call a supersubalgebra $\g\subset \osp(n|2k)$ {\it weakly-irreducible} if it does not preserve any non-degenerate
vector supersubspace of $\Real^n\oplus\Pi(\Real^{2k})$.

Let $\M$ and $\mathcal{N}$ be supermanifolds. Recall the definition of the product $\M\times \N=(M\times N,\O_{\M\times
\N})$. Let $U\subset M$ and $V\subset N$ be open subsets, and $(U,x^1,...,x^n,\xi^1,...,\xi^m)$ and
$(V,y^1,...,y^p,\eta^1,...,\eta^q)$ coordinate systems on $\M$ and $\N$, respectively. Then by definition, $\O_{\M\times
\N}(U\times V)=\O_{M\times N}(U\times V)\otimes\Lambda_{\xi^1,...,\xi^m,\eta^1,...,\eta^q}$ and this condition defines the
sheaf $\O_{\M\times \N}$ uniquely \cite{Leites,Leites1}.  Let $(\M,g)$ and $(\N,h)$ be Riemannian supermanifolds and let
$\na^g$ and $\na^h$ be the corresponding Levi-Civita connections. Then $g+h$ is a Riemannian supermetric on $\M\times\N$
and, obviously, $\hol(\M\times\N,g+h)_{(x,y)}=\hol(\M,g)_{x}\oplus\hol(\N,h)_{y}$ and
$\Hol(\t\na^{g+h})_{(x,y)}=\Hol(\t\na^g)_x\times \Hol(\t\na^h)_y$, where $x\in M$ and $y\in N$.  The following theorem
generalizes the Wu theorem \cite{Wu}.

\begin{theorem} Let $(\M,g)$ be a Riemannian supermanifold such that the pseudo-Riemannian manifold $(M,\tilde g)$ is simply connected
and geodesically complete. Then there exist Riemannian supermanifolds $(\M_0,g_0),(\M_1,g_1),...,(\M_r,g_r)$ such that
\begin{equation}\label{Mdec} (\M,g)=(\M_0\times\M_1\times\cdots\times\M_r,g_0+g_1+\cdots+g_r),\end{equation}
the supermanifold $(\M_0,g_0)$ is flat and the holonomy algebras of the supermanifolds $(\M_1,g_1)$,...,$(\M_r,g_r)$ are
weakly-irreducible. In particular,
\begin{equation}\hol(\M,g)=\hol(\M_1,g_1)\oplus\cdots\oplus\hol(\M_r,g_r).\end{equation}

For general $(\M,g)$ decomposition \eqref{Mdec} holds locally.
\end{theorem}

{\bf Proof.} The proof of the local version of this theorem is similar to the proof of the local version of the Wu
theorem. Let $x\in M$. Suppose that $\hol(\M,g)_x$ is not weakly-irreducible, then $\hol(\M,g)_x$ preserves a
non-degenerate vector supersubspace $F_1\subset T_x\M$. Let $F_2\subset T_x\M$ be its orthogonal complement. Then
$\hol(\M,g)_x$ preserves the decomposition $F_1\oplus F_2=T_x\M$. By Theorem \ref{paralsubsh}, there exist parallel
distributions $\F_1$ and $\F_2$ over on $\M$ defined over an open neighbourhood $U$ of the point $x$ (at the moment we do
not assume that $M$ is simply connected and that $\Hol(\t\nabla)_x$ preserves $F_1$ and $F_2$). As we have noted, the
distributions $\F_1$ and $\F_2$ are involutive. Hence, there exist maximal integral submanifolds $\M_1$ and $\M_2$ of $\M$
passing through the point $x$ and corresponding to the distributions $\F_1$ and $\F_2$, respectively \cite{Manin}.
Moreover, there exist local coordinates $x^1,...,x^{n},\xi^1,...,\xi^{m}$ (resp., $y^1,...,y^{n},\eta^1,...,\eta^{m}$) on
$\M$ such that $x^1,...,x^{n_1},\xi^1,...,\xi^{m_1}$ (resp., $y^1,...,y^{n-n_1},\eta^1,...,\eta^{m-m_1}$) are coordinates
on $\M_1$ (resp., on $\M_2$). Consequently,
$x^1,...,x^{n_1},y^1,...,y^{n-n_1},\xi^1,...,\xi^{m_1},\eta^1,...,\eta^{m-m_1}$ are coordinates on $\M$ and we see that
$\M$ is locally isomorphic to a domain in the product $\M_1\times \M_2$. Since $F_1$ and $F_2$ are non-degenerate, the
restrictions $g_1$ and $g_2$ of $g$ to $\F_1$ and $\F_2$, respectively, are non-degenerate.  It is easy to check that
$g_1$ and $g_2$ do not depend on the coordinates $y^1,...,y^{n-n_1},\eta^1,...,\eta^{m-m_1}$ and
$x^1,...,x^{n_1},\xi^1,...,\xi^{m_1}$, respectively. Thus, $(\M_1,g_1)$ and $(\M_2,g_2)$ are Riemannian supermanifolds and
$g=g_1+g_2$. The local version of the theorem is proved.

Suppose that the pseudo-Riemannian manifold $(M,\tilde g)$ is simply connected and geodesically complete. Let $F_1$,
$F_2$, $\F_1$, $\F_2$ $\M_1$ and $\M_2$ be as above. Obviously, the vector subspaces $(F_0)_{\bar 0},(F_1)_{\bar 0}\subset
T_xM$ are non-degenerate and preserved by the holonomy group $\Hol(M,\t g)_x$ of the pseudo-Riemannian manifold $(M,\t
g)$. By the Wu theorem, $M$ is diffeomorphic to the product $M_1\times M_2$, where $M_1$ and $M_2$ are integral
submanifolds passing through the point $x$ and corresponding to the parallel distributions defined by the vector subspaces
$(F_0)_{\bar 0}\subset T_xM$ and $(F_1)_{\bar 0}\subset T_xM$, respectively. It is obvious that the underlying manifolds
of the supermanifolds $\M_1$ and $\M_2$ are $M_1$ and $M_2$, respectively. From the local part of the theorem it follows
that $\M=\M_1\times \M_2$ and $g=g_1+g_2$. The theorem is proved. $\Box$

\section{Examples of Berger superalgebras}\label{exBerger}

In this section we give examples of complex Berger superalgebras. We use results and denotations of
\cite{EE,SoS,Poletaeva1,Poletaeva2}.

 Let $\g_{-1}$ denote a complex vector superspace and  let  $\g_0\subset\gl(\g_{-1})$ be a supersubalgebra.  {\it The $k$-th prolongation} ($k\geq 1$) $\g_k$ of $\g_0$ is defined as for representations of usual
Lie algebras up to additional signs. Consider {\it the Cartan prolong} $\g_*=\g_*(\g_{-1},\g_0)=\oplus_{k\geq -1}\g_k$.
Note that $\g_*$ has a structure of Lie superalgebra. By analogy with \cite{Sch} we get the following exact sequence
\begin{equation}\label{posled1}0\longrightarrow\g_2\longrightarrow\g_{-1}^*\otimes \g_1\longrightarrow\R(\g_0)\longrightarrow
H^{2,2}_{\g_0}\longrightarrow 0,\end{equation} where $H^{2,2}_{\g_0}$ is the $(2,2)$-th Spencer cohomology group (note
that this group is denoted in \cite{Sch} by $H^{1,2}_{\g_0}$). The second map in the sequence is given by
\begin{equation}\label{RVposled1}R_{\phi\otimes\alpha}(x,y)=\phi(x)\alpha(y)-(-1)^{|x||y|}\phi(y)\alpha(x).\end{equation}

In \cite{EE,Poletaeva1,Poletaeva2} examples of  irreducible subalgebras $\g_0\subset\gl(\g_{-1})$ with $\g_1\neq 0$ are
given and for the most of them the groups $H^{2,2}_{\g_0}$ are computed.

\begin{tab}\label{tab3} Examples of irreducible subalgebras $\g_0\subset\gl(\g_{-1})$ with $\g_1\neq 0$ and $\g_2=0$  $$
\begin{array}{|c|c|c|c|c|}\hline \g_0&\g_0:\g_{-1}&\g_0:\g_{1}&\g_*& \text{restriction}\\\hline

\mathfrak{c}(\sl(n-p|q)\oplus \sl(p|m-q))&\id^*\otimes\id&\g_{-1}^*&\sl(n|m)&n\neq m, n-p+q\geq 2,\\&&&&m-q+p\geq
2\\\hline

\sl(n-p|q)\oplus \sl(p|n-q)&\id^*\otimes\id&\g_{-1}^*&\mathfrak{psl}(n|n)&n\geq 3, n-p+q\geq 2,\\&&&&m-q+p\geq 2\\\hline

\mathfrak{cosp}(n|2k)  & \id&\id&\osp(n+2|2k)& \\\hline

\mathfrak{gl}(l|k)&\Lambda^2\id &\g_{-1}^*& \mathfrak{osp}(2l|2k)& (l,k)\neq (3,0)\\\hline

 \mathfrak{ps}(\mathfrak{q}(p)\oplus\mathfrak{q}(n-p))&\id^*\otimes\id&\g_{-1}^*&\mathfrak{psq}(n)&n\geq 3,\, n-1\geq p \geq 1\\\hline

\mathfrak{sl}(p|n-p)&\Pi(S^2\id)&\Pi(\Lambda^2\id^*)&\mathfrak{spe}(n)&n\geq 3,\ n\geq p\geq 0\\\hline

%\mathfrak{sl}(p|n-p)&\Pi(\Lambda^2\id)&\mathfrak{spe}(n)&n\geq 3,\ n\geq p\geq 0\\\hline

%\mathfrak{gl}(p|n-p)&\Pi(S^2\id)&\Pi(\Lambda^2\id^*)&\mathfrak{pe}(n)&n\geq 3,\ n\geq p\geq 0\\\hline

 %\mathfrak{gl}(p|n-p)&\Pi(\Lambda^2\id)&\mathfrak{pe}(n)&n\geq 3,\ n\geq p\geq 0\\\hline

\gl(1|2)&V_{(1+\alpha)\varepsilon_1}&V_{-\alpha\varepsilon_1}&D(\alpha)&\\\hline

\gl(1|2)&V_{\frac{1+\alpha}{\alpha}\varepsilon_1}&V_{-\frac{1}{\alpha}\varepsilon_1}&D(\alpha)&\\\hline

\gl(1|2)&V_{\frac{\alpha}{1+\alpha}\varepsilon_1}&V_{\frac{1}{1+\alpha}\varepsilon_1}&D(\alpha)&\\\hline

\mathfrak{cosp}(2|4) & V_{-\varepsilon_1+\delta_1+\delta_2}&V_{3\varepsilon_1}& \mathfrak{ab}(3)&\\\hline

%\mathfrak{vect}(0|n)\ltimes\Lambda(n)&\Pi(\Lambda(n))&\Pi(\mathfrak{vect}(0|n))&\mathfrak{vect}(0|n+1)&n\geq 1\\\hline

\end{array}$$\end{tab}

\begin{tab}\label{tab4} Examples of irreducible subalgebras $\g_0\subset\gl(\g_{-1})$ with $\g_1\neq 0$ and $\g_2\neq 0$
$$\begin{array}{|c|c|c|c|c|c|}\hline \g_0&\g_0:\g_{-1}&\g_0:\g_{1}&\g_0:\g_{2}&\g_*& \text{restriction}\\\hline

\mathfrak{cpe}(n)&\id&\id^*&\Pi(\Co)&\pe(n+1)&n\geq 3\\\hline

\mathfrak{spe}(n)\rtimes\langle \tau+nz\rangle &\id&\id^*&\Pi(\Co)&\mathfrak{spe}(n+1)&n\geq 3\\\hline

\end{array}$$\end{tab}

\begin{tab}\label{tab5} Examples of irreducible subalgebras $\g_0\subset\gl(\g_{-1})$ whose Caratan prolongs are vectorial Lie
superalgebras
  $$
\begin{array}{|c|c|c|c|c|}\hline \g_0&\g_0:\g_{-1}&\g_*&H_{\g_0}^{2,2}&\text{restriction}\\\hline
\gl(n|m)&\id&\mathfrak{vect}(n|m)&0&\\\hline \sl(n|m)&\id&\mathfrak{svect}(n|m)&0&(n,m)\neq (0,2)\\\hline
\sl(0|2)&\id&\mathfrak{svect}(0|2)&\Pi(\Co)&\\\hline \osp^{sk}(2n|m)&\id&\mathfrak{h}(2n|m)&0&\\\hline
\mathfrak{spe}^{sk}(n)&\id&\mathfrak{sle}(n)&\Pi(\Co)&n\geq 3\\\hline
\end{array}$$\end{tab}

\begin{prop}
The following Lie superalgebras are Berger superalgebras:
\begin{description} \item[1)] $\mathfrak{c}(\sl(n-p|q)\oplus \sl(p|m-q))$ and $\sl(n-p|q)\oplus \sl(p|m-q)$ if $n\neq m, n-p+q\geq 2,m-q+p\geq
2$, $\sl(n-p|q)\oplus \sl(p|n-q)$ if $n\geq 3,n-p+q\geq 2,n-q+p\geq 2$, $\mathfrak{cosp}(n|2k)$, $\mathfrak{osp}(n|2k)$,
$\mathfrak{ps}(\mathfrak{q}(p)\oplus\mathfrak{q}(n-p))$ and $\mathfrak{p}(\mathfrak{sq}(p)\oplus\mathfrak{sq}(n-p))$
 with their standard
representations;

\item[2)] $\mathfrak{gl}(l|k)$ and $\mathfrak{sl}(l|k)$ acting on  $\Lambda^2(\Real^l\oplus\Pi(\Real^k))$;

\item[3)] $\mathfrak{sl}(p|n-p)$  acting on both   $\Pi(S^2(\Real^p\oplus\Pi(\Real^{n-p})))$ and
$\Pi(\Lambda^2(\Real^p\oplus\Pi(\Real^{n-p})))$;

\item[4)] $\gl(1|2)$ and $\sl(1|2)$ acting on each of $V_{(1+\alpha)\varepsilon_1}$, $V_{-\alpha\varepsilon_1}$,
$V_{\frac{1+\alpha}{\alpha}\varepsilon_1}$, $V_{-\frac{1}{\alpha}\varepsilon_1}$
$V_{\frac{\alpha}{1+\alpha}\varepsilon_1}$ and $V_{\frac{1}{1+\alpha}\varepsilon_1}$;

\item[5)] $\mathfrak{cosp}(2|4)$ and $\mathfrak{osp}(2|4)$ acting on both $V_{-\varepsilon_1+\delta_1+\delta_2}$ and $V_{3\varepsilon_1}$;

%\item[6)] $\mathfrak{vect}(0|n)$ acting on $\Pi(\Lambda(n))$.

\end{description}

\end{prop}
{\bf Proof.} For the proof we use Table \ref{tab3}. In all cases $\g_*$ is simple. Hence, $[\g_{-1},\g_{1}]=\g_0$. This
and the exact sequence \eqref{posled1} yield that $\g_0$ is a Berger superalgebra: take $x,y\in \g_{-1}$, $\alpha\in\g_1$
and $\phi\in\g_{-1}^*$ such that $\alpha(x)\neq 0$, $\phi(y)=1$ and $\phi(x)=0$, then
$R_{\phi\otimes\alpha}(x,y)=\alpha(x)$; on the other hand, these elements span $\g_0$. Suppose that $\g_0$ is not
semisimple and $\g_0\neq\mathfrak{ps}(\mathfrak{q}(p)\oplus\mathfrak{q}(n-p))$, then   $\g_0=\hat\g_0\oplus\Co$ for an
ideal $\hat\g_0\subset\g_0$.   Let $\alpha\in (\g_{1})_{\bar 0}$. Then there exists  a non-zero $x\in(\g_{-1})_{\bar 0}$
such that $(\g_{-1})_{\bar 0}=\Co x\oplus\ker\pr_{\Co}\alpha|_{(\g_{-1})_{\bar 0}}$. Take a non-zero $\phi\in
(\g_{-1})_{\bar 0}$ such that $\phi(x)=0$. Than $0\neq R_{\phi\otimes\alpha}\in\R(\hat\g_0)$. Consequently,  if $\hat\g_0$
is simple then
 it is a Berger superalgebra. If $\hat\g_0$ is not simple, then it is clear that all elements of $\R(\hat\g_0)$ can not
take image in one of the simple summands of $\hat\g_0$. The case
$\g_0=\mathfrak{ps}(\mathfrak{q}(p)\oplus\mathfrak{q}(n-p))$ is similar.  $\Box$

\begin{prop} Let $n\geq 3$. The Lie subalgebras
$\mathfrak{spe}(n),\mathfrak{pe}(n),\mathfrak{cspe}(n),\mathfrak{cpe}(n),\mathfrak{spe}(n)\rtimes\langle\tau+ n
z\rangle,\mathfrak{spe}(n)\rtimes\langle a\tau+ b z\rangle$ ($a,b\in\Co,\,\frac{b}{a}\neq n$) of $\gl(n|n)$ are Berger
superalgebras.\end{prop}

{\bf Proof.} If $\g_0=\mathfrak{spe}(n),\mathfrak{pe}(n),\mathfrak{cspe}(n)$ or $\mathfrak{spe}(n)\rtimes\langle a\tau+ b
z\rangle$, then $\g_*(\g_{-1},\g_0)=0$. Moreover, $\R(\mathfrak{spe}(n))\simeq H^{2,2}_{\mathfrak{spe}(n)}$ and there are
the following non-split exact sequences $$0\longrightarrow V_{\varepsilon_1+\varepsilon_2} \longrightarrow
H^{2,2}_{\mathfrak{spe}(n)}\longrightarrow \Pi(V_{2\varepsilon_1+2\varepsilon_2})\longrightarrow 0,$$ where $n\geq 4$, and
$$0\longrightarrow V \longrightarrow H^{2,2}_{\mathfrak{spe}(3)}\longrightarrow \Pi(V_{3\varepsilon_1})\longrightarrow
0,$$ where $V$ is determined from the following non-split exact sequence $$0\longrightarrow
V_{\varepsilon_1+\varepsilon_2} \longrightarrow V\longrightarrow \Pi(V_{2\varepsilon_1+2\varepsilon_2})\longrightarrow
0.$$ Since $\mathfrak{spe}(n)$ is simple, it is a Berger superalgebra.
 For each of the last three values of $\g_0$ there is an
injective non-surjective  map from $H^{2,2}_{\mathfrak{spe}(n)}$ to $H^{2,2}_{\g_0}$.  Since each of these three Lie
superalgebras has dimension $\dim\mathfrak{spe}(n)+1$, they are Berger superalgebras.

Further, if $n\geq 4$, then $H^{2,2}_{\mathfrak{spe}(n)\rtimes\langle\tau+ n z\rangle}\simeq
\Pi(V_{2\varepsilon_1+2\varepsilon_2})$ and there is the following non-split exact sequence $$0\longrightarrow
\Pi(V_{2\varepsilon_1+2\varepsilon_2}) \longrightarrow H^{2,2}_{\mathfrak{spe}(3)\rtimes\langle\tau+ 3
z\rangle}\longrightarrow \Pi(V_{3\varepsilon_1})\longrightarrow 0.$$ From this, \eqref{posled1} and the exact sequence
$$0\longrightarrow\Pi(\Co)\longrightarrow S^2 \g_{-1}^*\longrightarrow V_{\varepsilon_1+\varepsilon_2}\longrightarrow 0$$
it follows that $\R(\mathfrak{spe}(n)\rtimes\langle\tau+ n z\rangle)\neq \R(\mathfrak{spe}(n))$. Hence
$\mathfrak{spe}(n)\rtimes\langle\tau+ n z\rangle$ is a Berger superalgebra. Finally, there is an injective non-surjective
map from $H^{2,2}_{\mathfrak{spe}(n)\rtimes\langle\tau+ n z\rangle}$ to $H^{2,2}_{\mathfrak{cpe}(n)}$. This and the
sequences \eqref{posled1} written for the both Lie superalgebras show that $\R(\mathfrak{spe}(n)\rtimes\langle\tau+ n
z\rangle)\neq\R(\mathfrak{cpe}(n))$. Since $\dim\mathfrak{cpe}(n)=\dim \mathfrak{spe}(n)\rtimes\langle\tau+ n z\rangle+1$,
we get that $\mathfrak{cpe}(n)$ is a Berger superalgebra. $\Box$

\begin{prop} The Lie subalgebras $\gl(n|m)$, $\sl(n|m)$, $\osp^{sk}(n|2m)$ and $\mathfrak{spe}^{sk}(k)$ ($k\geq 3$)
with their standard representations are Berger superalgebras.\end{prop}

{\bf Proof.} The proof follows from \eqref{posled1} and Table \ref{tab5}. $\Box$

\begin{prop}\label{Piadjoint} Let $\g_0$ be a simple complex
Lie superalgebra, $\g_{-1}=\Pi(\g_0)$ and $\g_0$ act  on $\g_{-1}$ vie the adjoint representation. Then $\g_1=\Co
\varphi_1$ and $\g_2=0$, where $\varphi_1:\Pi(\g_0)=\Pi((\g_0)_{\bar 1})\oplus\Pi((\g_0)_{\bar 0})\to\g_0=(\g_0)_{\bar
0}\oplus(\g_0)_{\bar 1}$, $\varphi_1(x)=(-1)^{|x|}\Pi(x)$ for all homogeneous $x\in\Pi(\g_0)$. In particular, for any
simple complex Lie superalgebra $\g$, the subalgebra $\g\subset\gl(\Pi(\g))$ is a Berger superalgebra.
\end{prop}

{\bf Proof.} We have \begin{multline}\g_1=\{\varphi\in
\Pi(\g_0)^*\otimes\g_0|\,\,[\varphi(x),y]=(-1)^{|x||y|}[\varphi(y),x],\,\, x,y\in\Pi(\g_0)\}\\=\Pi\{\varphi\in
\g_0^*\otimes\g_0|\,\,[\varphi(x),y]=(-1)^{(|x|+1)(|y|+1)}[\varphi(y),x],\,\, x,y\in\g_0\}.\end{multline} Let us first
find all $\varphi\in\g_1$ annihilated by $\g_0$. Suppose that $0\neq\varphi\in(\g_1)_{\bar 0}$ is annihilated by $\g_0$,
i.e. $[x,\varphi(y)]=\varphi([x,y])$ for all $x,y\in\g_0$. Consequently, the kernel of $\varphi$ is an ideal in $\g_0$ and
it must be trivial, i.e. $\varphi$ is injective. On the other hand, for $x,y\in\g_0$ we have
$\varphi([\varphi(x),y])=[\varphi(x),\varphi(y)]=-(-1)^{(|x|+1)(|y|+1)}[\varphi(y),\varphi(x)]=-(-1)^{(|x|+1)(|y|+1)}\varphi([\varphi(y),x])$,
since $\varphi\in(\Pi(\g_0)^*\otimes\g_0)_{\bar 0}=(\g_0^*\otimes\g_0)_{\bar 1}$. Hence, $\varphi([\varphi(x),y])=0$.
Since $\varphi$ is injective, $[\varphi(x),y]=0$. But this yields $\varphi=0$ and we get a contradiction. Let $\varphi_1$
be as in the statement of the proposition. Obviously, $\varphi_1\in(\g_1)_{\bar 1}$ and $\varphi_1$ is annihilated by
$\g_0$. From the above and the Schur Lemma it follows that the subset of $\g_1$ annihilated by $\g_0$ coincides with
$\Co\varphi_1$.

If $\g_1$ is not equal to $\Co\varphi_1$, then there exists a non-trivial $\g_0$-irreducible submodule $W\subset \g_1$.
Consider the $\mathbb{Z}$-graded Lie superalgebra $\h=\g_{-1}\oplus\g_0\oplus W\oplus W^2\oplus W^3\oplus\cdots$. By the
construction of $\h$ and a proposition from \cite{Kac}, $\h$ is simple. On the other hand, all simple $\mathbb{Z}$-graded
Lie superalgebras of depth 1 are listed in \cite{EE} and the case $\g_{-1}=\Pi(\g_0)$ does not occur there. Thus,
$\g_1=\Co\varphi_1$. $\Box$

{\bf Remark.} Note that Proposition \ref{Piadjoint} shows in particular that any simple vectorial Lie superalgebra $\g$
acting on $\Pi(\g)$ is a Berger superalgebra. These examples have no analogs in the case of the usual Berger algebras
\cite{Sch}.

\bibliographystyle{unsrt}

\vskip1cm

Department of Algebra and Geometry, Masaryk University in Brno, Jan\'a\v ckovo n\'am. 2a, 66295 Brno, Czech Republic

{\it E-mail address}: galaev@math.muni.cz

\end{document}